\documentclass[12pt]{article}

\usepackage{
	enumitem,
	times,
	amsmath,
	amssymb,
	amsthm,
	mathtools,
	fancyhdr,
	color,
	thmtools,
	etoolbox,
	xspace,
	hyperref,
	lscape,
	thm-restate
}

\usepackage[margin=2cm]{geometry}

\declaretheorem[style=definition,numberwithin=section]{definition}
\declaretheorem[style=definition,qed=$\oslash$,sibling=definition]{example}

\declaretheorem[style=slanted,sibling=definition]{theorem}
\declaretheorem[style=slanted,sibling=definition]{conjecture}
\declaretheorem[style=slanted,sibling=definition]{lemma}
\declaretheorem[style=slanted,sibling=definition]{proposition}
\declaretheorem[style=slanted,sibling=definition]{corollary}

\declaretheorem[style=definition,sibling=example]{remark}

\AtEndEnvironment{proof}{\setcounter{claim}{0}}

\usepackage[square, comma, numbers, sort&compress ]{natbib}
\newcommand{\C}{\ensuremath{\mathbb{C} }}
\newcommand{\R}{\ensuremath{\mathbb{R} }}
\newcommand{\N}{\ensuremath{\mathbb{N} }}
\newcommand{\Z}{\ensuremath{\mathbb{Z} }}
\DeclareMathOperator{\Tr}{Tr}
\DeclareMathOperator{\End}{End}
\DeclareMathOperator{\Ad}{Ad}
\DeclareMathOperator{\ad}{ad}

\newcommand\scalemath[2]{\scalebox{#1}{\mbox{\ensuremath{\displaystyle #2}}}}
\title{On the Mathieu Conjecture for $Sp(N)$ and $G_2$}
\author{Kevin Zwart\footnote{Email address: \href{mailto:kevin.zwart@ru.nl}{kevin.zwart@ru.nl}}}
\date{%
	Radboud University Nijmegen, IMAPP \\[2ex]%
	}

\begin{document}

\maketitle
\abstract{As a direct continuation of K. Zwart \cite{Zwart}, which is built on the work of M. M\"uger and L. Tuset, we reduce the Mathieu conjecture, formulated by O.\ Mathieu in 1997, for $Sp(N)$ and $G_2$ to a conjecture involving functions over $\R^n\times (S^1)^m$ with $n,m\in\N_0$. The proofs rely on Euler-style parametrizations of these groups, a specific version of the $KAK$ decomposition, which we discuss and prove.}      

\section{Introduction}

In 1997, Olivier Mathieu conjectured the following statement:
\begin{conjecture}\cite{Mathieu}
	Let $G$ be a compact connected Lie group. If $f,h$ are finite-type functions such that $\int_G f^P \,dg=0$ for all $P\in\N$, then $\int_G f^Ph \,dg=0$ for all large enough $P$.
\end{conjecture}

Just one year after the publication of Mathieu's paper, Duistermaat and van der Kallen \cite{Duistermaat} proved Mathieu's conjecture for all abelian connected compact groups. Although the conjecture is still an open problem for any non-abelian group, some attempts were made. Dings and Koelink \cite{Dings-Koelink} approached the conjecture for $SU(2)$ by expressing the finite-type functions in terms of explicit matrix coefficients. Influenced heavily by this, M\"uger and Tuset \cite{Mueger} reduced the Mathieu conjecture for $SU(2)$ to a conjecture about certain Laurent polynomials. In 2023, this author proved \cite{Zwart} that a similar reduction can be done for $SU(N)$ and $SO(N)$ to a conjecture regarding Laurent polynomials which also allow for complex square roots.

This paper is a direct continuation of \cite{Zwart} where we apply the same arguments to the infinite family of groups $Sp(N)$ with $N\in\N$ and the exceptional group $G_2$. As in \cite{Zwart}, this can be achieved by using a generalization of the Euler decomposition which we will call the Euler angles. This decomposition is a more intricate version of the $KAK$ decomposition, a well-known decomposition for connected Lie groups. The decomposition will be proven and applied to the groups $Sp(N)$ and $G_2$.

After having found the Euler angles, the finite-type functions of $Sp(N)$ and $G_2$ are considered, and we will reduce the Mathieu conjecture to a similar conjecture to that of M\"uger and Tuset in \cite{Mueger} and Zwart \cite{Zwart} with a different weight function. To be more precise, any finite type function of $Sp(N)$ or $G_2$ reduces to a function $f:\C^n\times \R^k\rightarrow\C$ which can be written as $f(z,x)=\sum_{\vec{m}}c_{\vec{m}}(x)z^{\vec{m}}$ where $\vec{m}=(m_1,\ldots,m_k)$ is a multi-index where $m_i\in\bigcup_{j=1}^l\frac{1}{j}\mathbb{Z}$ for each $i$, for some $l\in\N$, and $c_{\vec{m}}(x)$ is a polynomial in $x_1,\ldots,x_k$ and $\sqrt{1-x_1^2},\ldots,\sqrt{1-x_k^2}$. Assuming these functions satisfy certain conjectures, the Mathieu conjecture is proven for $Sp(N)$ and $G_2$. In Section \ref{sec:Sp(N)} we will focus on the group $Sp(N)$, while in Section \ref{sec:G2} the group $G_2$ will be considered. The final part of the paper is dedicated to proving the generalized Euler angles we use throughout this paper, with the corresponding explicit description of the Haar measure in this parametrization.

\textit{Acknowledgment:} The author would like to thank Michael M\"uger for the continued help with the project, and the valuable discussions. He also wishes to thank Erik Koelink for the suggestion of looking at these two specific groups, and the valuable support.

\section{The case of $Sp(N)$}\label{sec:Sp(N)}
\subsection{The Euler angles decomposition of $Sp(N)$}
In this paper we will discuss Mathieu's conjecture on $Sp(N)$ and $G_2$ with $N \in \N$. We start by recalling Mathieu's conjecture. To do so, we first introduce the notion of a finite-type function:
\begin{definition}
	Let $G$ be a compact Lie group. A function $f:G\rightarrow\C$ is called a \emph{finite-type function} if it can be written as a finite sum of matrix components of irreducible continuous representations, i.e. $$f(x)=\sum_{j=1}^N\Tr\left[a_j\pi_j(x)\right],$$ where $(\pi_j,V_j)$ is an irreducible continuous representation of $G$, and  $a_j\in\End(V_j)$.  
\end{definition}

\begin{conjecture}[The Mathieu Conjecture \cite{Mathieu}]
	Let $G$ be a compact connected Lie group. If $f,h$ are finite-type functions such that $\int_G f^P \,dg=0$ for all $P\in\N$, then $\int_G f^Ph \,dg=0$ for all large enough $P$.
\end{conjecture}

This section is dedicated to the Lie group $Sp(N)$, while the next section is dedicated to $G_2$. This document is a direct continuation of \cite{Zwart}, and we assume the reader is familiar with the techniques used there, for we will apply them throughout the paper. 

In a similar fashion as in \cite{Zwart}, we wish to apply some form of $KAK$ decomposition to our groups in such a way that we can reduce our Mathieu conjecture. As already hinted in \cite{Zwart}, the generalized Euler angles are an example of a more general decomposition, which is stated below in Theorem \ref{thm:KAK_decomp}. This theorem is key to the rest of the paper, and will be applied to both $Sp(N)$ and $G_2$. We will prove the theorem in Appendix \ref{appendix:Proof_General_KAK_decomp} to keep this section focussed on $Sp(N)$.

\begin{restatable}{theorem}{KAKDecomp}\emph{[The Euler Angles Theorem]}\label{thm:KAK_decomp}
	Let $G$ be a simply connected compact Lie group with finite center, and let $\mathfrak{g}$ be its Lie algebra. Let $\theta:\mathfrak{g}\rightarrow\mathfrak{g}$ be an involutive automorphism. Let $\mathfrak{k}, \mathfrak{p}$ be the $+1$ and $-1$ eigenspace of $\theta$, respectively, in such a way that $\mathfrak{g}=\mathfrak{k}\oplus\mathfrak{p}$. Fix a maximal abelian subalgebra $\mathfrak{a}\subset \mathfrak{p}$ and fix a set of (positive) roots $\Delta_\mathfrak{p}$ with respect to $\mathfrak{a}$. Then the mapping 
	\begin{align}
		(K/M)\times \exp(\mathcal{A})\times K \rightarrow G,\qquad\qquad (kM,\exp(H),l)\mapsto \exp(\Ad_G(k)H)l
	\end{align}
	is surjective and a diffeomorphism up to a measure zero set if we replace $\mathcal{A}$ by $\mathrm{int}(\mathcal{A})$. Here $K\subseteq G$ is the connected Lie group with Lie algebra $\mathfrak{k}$, and $M=Z_K(\mathfrak{a})$. Here $\mathcal{A}$ is the closure of the connected component of the following set
	\begin{align*}
		\mathfrak{a}_{+}-\{H\in\mathfrak{a}\,|\,\alpha(H)\in \pi i\Z \,\text{ for some } \alpha\in \Delta_\mathfrak{p} \}
	\end{align*} in such a way that $0\in \mathcal{A}$. Here $\mathfrak{a}_+$ is the closed positive Weyl chamber in $\mathfrak{a}$.
	In addition, the Haar measure decomposes in this parameterization as follows
	\begin{align*}
		\int_G h(g) dg = C \int_{K/M}\int_{\mathcal{A}}\int_K h(k_1\exp(H)k_2) |J(\exp(H))| dk_2 dH dg_{K/M}
	\end{align*}
	for any measurable function $h:G\rightarrow \C$, where $C>0$ is a constant (independent of $h$), $k_2\in K$ with corresponding Haar measure $dk_2$, $k_1\in k_1M$ an arbitrary representative of $k_1M\in K/M$ with corresponding unique $K$-invariant measure $dg_{K/M}$ on $K/M$, and $dH$ is the measure on $\mathfrak{a}$. Here $J:A\rightarrow\C$ is given by
	\begin{align}\label{eq:Jacobian_abstract_Euler_decomp}
		J(\exp(H)) := \prod_{\alpha\in \Delta_\mathfrak{p}^+} \sin(\alpha(iH))
	\end{align}
	where $\Delta_\mathfrak{p}^+$ is the set of positive roots on $\mathfrak{a}$.
\end{restatable}


The rest of the section will be dedicated to applying Theorem \ref{thm:KAK_decomp} to the group $Sp(N)$. To do that, we first recall the definition of the compact symplectic group $Sp(N)$.

\begin{definition}
	Let $N\in\N$. We define the \emph{compact symplectic group} as
	\begin{align*}
		Sp(N):=\{A\in U(2N) \,|\, JAJ^T=A \}
	\end{align*}
	where $J=\begin{pmatrix}
		0 & \mathbf{1}_N\\
		-\mathbf{1}_N &0
	\end{pmatrix}$ and $\mathbf{1}_N$ is the $N \times N$ identity matrix.
\end{definition}
It is clear that the Lie algebra of $Sp(N)$ is given by
\begin{align}
	\mathfrak{sp}(N)=\left\{\left.\begin{pmatrix}
		U & A\\
		-A^\dagger & U^*
	\end{pmatrix}\right| U\in \mathfrak{u}(N) \text{ and } A\in \mathfrak{gl}(N,\C) \text{ is a symmetric matrix}\right\}.
\end{align}
Here, we denoted $A^\dagger$ as the adjoint of $A$, and $U^*$ the complex conjugate of $U$. In addition, since $\mathfrak{sp}(N)\subset\mathfrak{u}(2N)$, the Killing form is just given by the bilinear form $\langle X,Y\rangle=4\Tr(XY)$.
We note that the map $$\theta:=\Ad\begin{pmatrix}
i\mathbf{1}_N&0\\
0&-i\mathbf{1}_N
\end{pmatrix}$$ has the property $\theta\in \mathrm{Int}(\mathfrak{sp}(N))$ and can be seen as a Cartan involution. It is easily seen that with this Cartan involution, we get the following decomposition:
\begin{align}
	\mathfrak{k}&:= \left\{\left.\begin{pmatrix}
	U&0\\
	0&U^*
	\end{pmatrix}\right| U\in \mathfrak{u}(N) \right\}\simeq \mathfrak{u}(N),\\
	\mathfrak{p}&:= \left\{\left.\begin{pmatrix}
	0&A\\
	-A^\dagger&0
	\end{pmatrix}\right| A\in \mathfrak{gl}(N,\C)\text{ such that } A \text{ is symmetric}\right\}.
\end{align} which gives $$K:=\left\{\left.\begin{pmatrix}
U&0\\
0&U^*
\end{pmatrix}\right| U\in U(N) \right\}\simeq U(N).$$
One can take the following set as maximal abelian subalgebra:
\begin{align*}
	\mathfrak{a}:=\left\{\left.\begin{pmatrix}
	\mathbf{0}_N&\begin{matrix}
	a_1&&&\\
	&a_2&&\\
	&&\ddots&\\
	&&&a_N
	\end{matrix}\\
	\begin{matrix}
	-a_1&&&\\
	&-a_2&&\\
	&&\ddots&\\
	&&&-a_N
	\end{matrix}&\mathbf{0}_N
	\end{pmatrix}\right| a_j\in \R\text{ for all }j \right\}
\end{align*}
where $\mathbf{0}_N$ is the $N\times N$-matrix consisting of only zeros. With that, we can calculate $M=Z_K(\mathfrak{a})$. Since $M\subseteq K$, and it must commute with all options of $a_j$, we see that the only options for $x\in M$ is to be of the form $$x=\mathrm{diag}(e^{i\phi_1},\ldots,e^{i\phi_N},e^{-i\phi_1},\ldots,e^{-i\phi_N})$$ for some $\phi_j\in\R$. In specific, we see that
\begin{align*}
	\Ad(x)\begin{pmatrix}
	\mathbf{0}_N&\begin{matrix}
	a_1&&\\
	&\ddots&\\
	&&a_N
	\end{matrix}\\
	\begin{matrix}
	-a_1&&\\
	&\ddots&\\
	&&-a_N
	\end{matrix}&\mathbf{0}_N
	\end{pmatrix} = \begin{pmatrix}
	\mathbf{0}_N&\begin{matrix}
	e^{2i\phi_1}a_1&&\\
	&\ddots&\\
	&&e^{2i\phi_N}a_N
	\end{matrix}\\
	\begin{matrix}
	-e^{-2i\phi_1}a_1&&\\
	&\ddots&\\
	&&-e^{-2i\phi_N}a_N
	\end{matrix}&\mathbf{0}_N
	\end{pmatrix}.
\end{align*}
Hence $\phi_j=0$ or $\pi$ for all $j$, giving 
\begin{align*}
	M&=\langle\mathbf{1}_N,\mathrm{diag}(-1,1,\ldots,1,-1,1,\ldots,1),\mathrm{diag}(1,-1,\ldots,1,1,-1,\ldots,1),\ldots,\mathrm{diag}(1,1,\ldots,-1,1,1,\ldots,-1)\rangle\\
	&\simeq \Z_2^N
\end{align*}
Therefore, applying Theorem \ref{thm:KAK_decomp}, we find that
\begin{align*}
	Sp(N)&=(K/M) A\, K\\
	&\simeq (U(N)/\mathbb{Z}_2^N)\, \exp(\mathcal{A})\, U(N).
\end{align*}
For the rest of this paper we will work with the latter diffeomorphism, since on $U(N)$ some results regarding the Euler angle decomposition are already known \cite{Zwart}. Let us consider the parametrization of $U(N)$ first.

We remind ourselves of the Euler angle decomposition of $SU(N)$ and on $U(N)$. To describe this decomposition, we require an orthogonal basis for $\mathfrak{u}(N)$. Let $j=1,2,\ldots,N-1$ and $k=1,2,\ldots,2j$ and define the matrices $\lambda_0,\lambda_j\in \mathfrak{su}(N)$ in the following way\footnote{In most physics papers the matrices $\{i\lambda_j\}_j$ with $j>0$ are called Gell-Mann matrices, see e.g. \cite{Parametrization_SU(N),Parametrization_SU(4),bertini2006euler}}
\begin{align*}
[\lambda_0]_{ij}&:=i\delta_{i,1}\delta_{j,1}\\
[\lambda_{j^2-1+k}]_{\mu,\nu}&:=i(\delta_{\lceil\frac{k}{2}\rceil,\mu}\delta_{j+1,\nu}+\delta_{j+1,\mu}\delta_{\lceil\frac{k}{2}\rceil,\nu})\quad\qquad\text{if $k$ is odd},\\	[\lambda_{j^2-1+k}]_{\mu,\nu}&:=\delta_{\frac{k}{2},\mu}\delta_{j+1,\nu}-\delta_{j+1,\mu}\delta_{\frac{k}{2},\nu}\qquad\qquad\quad\text{if $k$ is even},\\ 
[\lambda_{(j+1)^2-1}]&:=\begin{pmatrix}
i\mathbf{1}_j&0&&&\\
0&-ij&&&\\
&&0&&\\
&&&\ddots&\\
&&&&0
\end{pmatrix}.
\end{align*}
This is an orthonormal basis for $\mathfrak{u}(N)$ as can easily be checked. With that, we get the following lemma:
\begin{restatable}[Generalized Euler Angles $SU(N)$, \cite{Zwart}]{lemma}{EulerAngles}\label{lemma:Euler_Angles}
	Let $N\geq 1$. Define inductively the mapping $F_{SU(N)}:([0,\pi]\times[0,2\pi]^{N-2})\times([0,\pi]\times[0,2\pi]^{N-3})\times\cdots\times([0,\pi]\times [0,2\pi])\times[0,\pi]\times \left[0,\frac{\pi}{2}\right]^{\frac{N(N-1)}{2}}\times[0,2\pi]\times\cdots\times\left[0,\frac{2\pi}{N-1}\right]\rightarrow SU(N)$ by $F_{SU(1)}\equiv 1$ and
	\begin{equation}\label{eq:Euler_parametrization}
	\begin{split}
	&F_{SU(N)}(\phi_1,\ldots\phi_{\frac{N(N-1)}{2}},\psi_1,\ldots,\psi_{\frac{N(N-1)}{2}},\omega_1,\ldots,\omega_{N-1}):=\\
	&\left(\prod_{2\leq k\leq N}A(k)(\phi_{k-1},\psi_{k-1})\right)\cdot\begin{pmatrix}
	F_{SU(N-1)}(\phi_{N},\ldots,\phi_{\frac{N(N-1)}{2}},\psi_{N},\ldots,\psi_{\frac{N(N-1)}{2}},\omega_1,\ldots,\omega_{N-2})&0\\
	0&1
	\end{pmatrix} e^{\lambda_{N^2-1}\omega_{N-1}},
	\end{split}	
	\end{equation}
	where $A(k)(x,y):=e^{\lambda_{3}x}e^{\lambda_{(k-1)^2+1}y}$, and $\psi_j\in\left[0,\frac{\pi}{2}\right],\,\omega_j\in\left[0,\frac{2\pi}{j}\right]$ for all $j$. Here we denoted the product as $$\prod_{2\leq k\leq N}A(k)(\phi_{k-1},\psi_{k-1}):=A(2)(\phi_1,\psi_1)\cdot \cdots\cdot A(N)(\phi_{N-1},\psi_{N-1}).$$ This mapping is surjective. Moreover it is a diffeomorphism on the interior of the hypercube onto its image which is $SU(N)$ up to a measure zero set.
\end{restatable}

\begin{corollary}[Generalized Euler Angles $U(N)$]\label{cor:KAK_decomp_U(N)}
	Let $N\in\N$. Then the mapping $F_{U(N)}: ([0,\pi]\times[0,2\pi]^{N-2})\times([0,\pi]\times[0,2\pi]^{N-3})\times\cdots\times([0,\pi]\times [0,2\pi])\times[0,\pi]\times \left[0,\frac{\pi}{2}\right]^{\frac{N(N-1)}{2}}\times[0,2\pi]\times\cdots\times\left[0,\frac{2\pi}{N-1}\right]\times [0,2\pi]\rightarrow U(N)$ defined by $$F_{U(N)}(\phi_1,\ldots,\omega_{N-1},\xi):=F_{SU(N)}(\phi_1,\ldots,\omega_{N-1})e^{\xi\lambda_0}$$ is a surjective map onto $U(N)$, and is a diffeomorphism on the interior of the hypercube onto its image which is $U(N)$ up to a measure zero set.
\end{corollary}
\begin{proof}
	This can be seen by the fact that $SU(N)\times U(1)\simeq U(N)$ as manifolds by sending $(A,z)\mapsto A\cdot \mathrm{diag}(z,1,\ldots,1)$ with inverse $A\mapsto \left(A\cdot \mathrm{diag}(\det(A)^{-1},1,\ldots,1),\det(A)\right)$. Using that $x\mapsto e^{ix}$ with $x\in[0,2\pi]$ is a diffeomorphism on the interior of the hypercube, and $\mathrm{diag}(z,1,...,1)=e^{ix\lambda_0}$ in this parametrization, gives the corollary.
\end{proof}

Now that we know how to parametrize $U(N)$, we go to $U(N)/\mathbb{Z}_2^N$. To do so, note that the $N$ rightmost matrices in $F_{U(N)}(\phi_1,\ldots,\xi)$ are diagonal matrices. Since the symmetric space $U(N)/\mathbb{Z}_2^N$ is found by considering the coset elements $g\mathbb{Z}_2^N$ which is a finite set, and these sets are found by multiplying diagonal matrices on the right, we can argue that the $\omega_1,\ldots,\omega_{N-1},\xi$ parameters will have a shorter range which will fully describe $U(N)/\mathbb{Z}_2^N$. In other words we can see $U(N)/\mathbb{Z}_2^N$ as a submanifold of $U(N)$. To be more precise, we note that
\begin{align*}
	\left.e^{\xi\lambda_{0}}\right|_{\xi=\pi}=\begin{pmatrix}-1&&&\\
	&1&&\\
	&&\ddots&\\
	&&&1	
	\end{pmatrix}\in \mathbb{Z}_2^N, \qquad \left.e^{\omega_1\lambda_3}e^{\xi\lambda_0}\right|_{\omega_1=\pi,\xi=\pi}=\begin{pmatrix}
	1&&&\\
	&-1&&\\
	&&\ldots&\\
	&&&1
	\end{pmatrix}\in \mathbb{Z}_2^N
\end{align*}
and so on. We conclude that we can find a parametrization of all elements in $\mathbb{Z}_2^N$ this way, hence giving the following parametrization:

\begin{corollary}[Parametrization $U(N)/\mathbb{Z}_2^N$.]\label{cor:KAK_decomp_U(N)/Z_2^N}
	Let $N\in \N$. Then the mapping $\tilde{F}_{U(N)}: ([0,\pi]\times[0,2\pi]^{N-2})\times([0,\pi]\times[0,2\pi]^{N-3})\times\cdots\times([0,\pi]\times [0,2\pi])\times[0,\pi]\times \left[0,\frac{\pi}{2}\right]^{\frac{N(N-1)}{2}}\times[0,\pi]\times\cdots\times\left[0,\frac{\pi}{N-1}\right]\times [0,\pi]\rightarrow U(N)$ defined by $$\tilde{F}_{U(N)}(\phi_1,\ldots,\omega_{N-1},\xi):=F_{SU(N)}(\phi_1,\ldots,\omega_{N-1})e^{\xi\lambda_0}$$ is a surjective map onto $U(N)/\mathbb{Z}_2^N$, and is a diffeomorphism on the interior of the hypercube onto its image which is $U(N)/\mathbb{Z}_2^N$ up to a measure zero set.
\end{corollary}
With this corollary, we get the following lemma which parametrizes $Sp(N)$. We will call this the \emph{Euler angles parametrization} or the \emph{Euler angles decomposition} of $Sp(N)$

\begin{lemma}[Euler Angles Parametrization of $Sp(N)$]\label{lemma:KAK_decomp_Sp(N)}
	Let $N\in\N$ and define $G=Sp(N)$, fix a maximal abelian subalgebra $\mathfrak{a}\subseteq\mathfrak{p}$ and consider $K\simeq U(N)$ as before. Then $$G\simeq (K/M)\,A\,K$$ diffeomorphic up to a measure zero set, where $A=\exp(\mathcal{A})\subseteq \exp(\mathfrak{a})$ and $M=Z_K(\mathfrak{a})$. Here 
	\begin{align*}
	\mathcal{A}=\left\{(y_1,\ldots,y_N)\in\R^N\left|\,0\leq y_j\leq \frac{\pi}{2}\text{ and } y_i\leq y_{i+1} \text{ for all }i=1,\ldots N-1 \right\}\right..
	\end{align*}
	In addition, let $F_{U(N)}$ and $\tilde{F}_{U(N)}$ be as in Corollary \ref{cor:KAK_decomp_U(N)} and Corollary \ref{cor:KAK_decomp_U(N)/Z_2^N} respectively. Then we define the mapping $F_{Sp(N)}:([0,\pi]\times[0,2\pi]^{N-2})\times([0,\pi]\times[0,2\pi]^{N-3})\times\cdots\times([0,\pi]\times [0,2\pi])\times[0,\pi]\times \left[0,\frac{\pi}{2}\right]^{\frac{N(N-1)}{2}}\times[0,\pi]\times\cdots\times\left[0,\frac{\pi}{N-1}\right]\times [0,\pi] \times \mathcal{A} \times ([0,\pi]\times[0,2\pi]^{N-2})\times([0,\pi]\times[0,2\pi]^{N-3})\times\cdots\times([0,\pi]\times [0,2\pi])\times[0,\pi]\times \left[0,\frac{\pi}{2}\right]^{\frac{N(N-1)}{2}}\times[0,2\pi]\times\cdots\times\left[0,\frac{2\pi}{N-1}\right]\times [0,2\pi] \rightarrow Sp(N)$ by
	\begin{align*}
		&F_{Sp(N)}(\tilde{\phi}_1,\ldots,\tilde{\omega}_{N-1},\tilde{\xi},y_1,\ldots,y_N,\phi_1,\ldots\omega_{N-1},\xi):=\\
		&\begin{pmatrix}
		\tilde{F}_{U(N)}(\tilde{\phi}_1,\ldots,\tilde{\omega}_{N-1},\tilde{\xi})&0\\
		0&\tilde{F}_{U(N)}(\tilde{\phi}_1,\ldots,\tilde{\omega}_{N-1},\tilde{\xi})^*
		\end{pmatrix} \exp\begin{pmatrix}
		\mathbf{0}_N&\begin{matrix}
		y_1&&\\
		&\ddots&\\
		&&y_n
		\end{matrix}\\
		\begin{matrix}
		-y_1&&\\
		&\ddots&\\
		&&-y_n
		\end{matrix}&\mathbf{0}_N
		\end{pmatrix}\cdot\\
		&\qquad\qquad\qquad\qquad\qquad\qquad\qquad\qquad\qquad\qquad\qquad \begin{pmatrix}
		F_{U(N)}(\phi_1,\ldots,\omega_{N-1},\xi)&0\\
		0&F_{U(N)}(\phi_1,\ldots,\omega_{N-1},\xi)^*
		\end{pmatrix}
	\end{align*}
	This mapping is surjective onto $Sp(N)$ and is a diffeomorphism on the interior of the hypercube onto its image which is $Sp(N)$ up to a measure zero set.
	
	Finally, the Haar measure $dg$ is in this parametrization given by:
	\begin{align*}
		dg&= dg_{K/M}(\tilde{\phi}_1,\ldots,\tilde{\xi})dg_{K}(\phi_1,\ldots,\xi)J(y_1,\ldots,y_N) dy_1\ldots dy_N\\
		& = C dg_{K}(\tilde{\phi}_1,\ldots,\tilde{\xi})dg_{K}(\phi_1,\ldots,\xi)J(y_1,\ldots,y_N) dy_1\ldots dy_N
	\end{align*} where $C>0$ is some constant, $dg_{K}(\phi_1,\ldots,\xi)$ is the Haar measure on $K\simeq U(N)$ integrating over the variables $\phi_1,\ldots,\xi$, and 
	\begin{align*}
		J_{Sp(N)}(y_1,\ldots,y_N):=\prod_{j=1}^{N}\sin(2y_j)\prod_{i> k}\sin(y_i-y_k)\sin(y_i+y_k).
	\end{align*}
	Here the integration over the $y_j$ variables goes over exactly the area $\mathcal{A}$.
\end{lemma}

\begin{remark}
	Note that we chose the variables that go over the symmetric space $K/M$ to have a $\sim$ above their letters for bookkeeping, while the variables going over elements in $K$ or $A$ are without a $\sim$. This notation is kept throughout the paper.
\end{remark}

\begin{proof}
	The first equation is proven by applying Theorem \ref{thm:KAK_decomp} to our case of $Sp(N)$. Secondly, by applying Corollary \ref{cor:KAK_decomp_U(N)} and Corollary \ref{cor:KAK_decomp_U(N)/Z_2^N} to the first equation, we find that $F_{Sp(N)}(\tilde{\phi},\ldots,\xi)$ is a surjective map onto $Sp(N)$ and is a diffeomorphism up to a measure zero set.
	The only thing to prove is the area $\mathcal{A}$. To do so, we recall the reduced root system of $\mathfrak{sp}(N)$ with respect to $\mathfrak{a}$. Let us define the basis of $\mathfrak{a}$ as the set $\{e_j\}_{j=1,\ldots,N}$ where
	\begin{align*}
		[e_j]_{k,l}=\delta_{k,j}\delta_{l,N+j}-\delta_{l,j}\delta_{k,N+j}.
	\end{align*}
	Define the linear functional $\alpha_j\in\mathfrak{a}^*$ as
	\begin{align*}
		\alpha_j(e_k)=\delta_{j,k}
	\end{align*}
	Then $\mathfrak{sp}(N)_\C$ has a root decomposition, with roots given by $$\Delta_{\mathfrak{p}}=\{\pm 2\alpha_j,\alpha_j\pm \alpha_k,-\alpha_j\pm \alpha_k\,|\, j,k=1,\ldots,N, j > k\}.$$ We choose the set of positive roots to be $$\Delta_\mathfrak{p}^+:=\{2\alpha_j,\alpha_j\pm \alpha_k\,|\,j>k\}.$$ Then it is immediate that 
	\begin{align*}
		J_{Sp(N)}(\exp(H))&=\prod_{j=1}^N\sin(2\alpha_j(H))\prod_{j>k}\sin((\alpha_j-\alpha_k)(H))\sin((\alpha_j+\alpha_k)(H))\\
		&=\prod_{j=1}^N\sin(2y_j)\prod_{j>k}\sin(y_j-y_k)\sin(y_j+y_k)
	\end{align*}
	where we wrote 
	\begin{align*}
		H=\begin{pmatrix}
		\mathbf{0}_N&\begin{matrix}
		y_1&&\\
		&\ddots&\\
		&&y_n
		\end{matrix}\\
		\begin{matrix}
		-y_1&&\\
		&\ddots&\\
		&&-y_n
		\end{matrix}&\mathbf{0}_N
		\end{pmatrix}.
	\end{align*} Thus $\mathcal{A}$ can be described as
	\begin{align*}
		\mathcal{A}=\left.\left\{(y_1,\ldots,y_N)\in\R^N\right|0\leq 2y_j\leq \pi\text{ and } 0\leq y_j\pm y_k\leq \pi \text{ for all }j>k \right\}.
	\end{align*}
	Note that if $0\leq y_j\leq\frac{\pi}{2}$ then $y_j\pm y_k\leq\pi$ automatically, so we can drop most inequalities, and only $0\leq y_j-y_k$ adds new information. In other words, we can describe $\mathcal{A}$ as
	\begin{align*}
		\mathcal{A}&=\left\{(y_1,\ldots,y_N)\in\R^N\left|\,0\leq 2y_j\leq \pi\text{ and } y_k\leq y_j \text{ for all }j>k \right\}\right.\\
		&=\left\{(y_1,\ldots,y_N)\in\R^N\left|\,0\leq 2y_j\leq \pi\text{ and } y_i\leq y_{i+1} \text{ for all }i=1,\ldots N-1 \right\}\right..
	\end{align*}
	 which yields the theorem.
We finally note that the measure on $K/M$ can be seen to be the Haar measure on $K$ times some positive constant. This can be seen by noting that up to a measure zero set, we can embed $K/M$ into $K$. This allows us to note that the left-invariant measure on $K/M$ can then be seen as the measure on $K$ as well up to a $C^\infty$-function. But since both $dg_{K/M}$ and $dg_K$ are left-invariant, this $C^\infty$ function should be constant. 
\end{proof}

\subsection{The Mathieu Conjecture on $Sp(N)$}\label{sec:Mathieu_Sp(N)}

Now that we have a suitable parametrization for $Sp(N)$, we go to the Mathieu conjecture. We work in the same way as in \cite{Zwart}. By Procesi \cite[Thm. 8.2.3]{Procesi}, all finite-type functions of compact groups $G\subset GL(N,\C)$ are generated by the matrix entries of $G$ and the inverse of the determinant. With our parametrization of $Sp(N)$ in Theorem \ref{lemma:KAK_decomp_Sp(N)}, we can thus describe all finite-type functions almost everywhere. Together with the fact that the determinant of all matrices in $Sp(N)$ is one, we see that any finite type function on $Sp(N)$ can be written as a sum of products of the matrix entries of $Sp(N)$, which is by Theorem \ref{lemma:KAK_decomp_Sp(N)} is described by the matrix entries of $K/M$, $\exp\mathcal{A}$ and $K$. Since the matrix entries of $K$ (and of $K/M$) are by Corollary \ref{cor:KAK_decomp_U(N)} given by the matrix entries of $SU(N)$ times the determinant of $U(N)$, which again can be interpreted as finite-type functions on $SU(N)$ times the determinant of $U(N)$, we get that a finite type function $f$ on $Sp(N)$ are can be written as

\begin{align}\label{eq:Finite_type_Sp(N)}
	f(g)=\sum_{i,j} c_{ij} f^{SU(N)}_{ij}(\tilde{\phi}_1,\ldots,\tilde{\omega}_{N-1})e^{im_{ij}\tilde{\xi}}\,g_{ij}(y_1,\ldots,y_N)\,h^{SU(N)}_{ij}(\phi_1\ldots,\xi)e^{in_{ij}\xi}
\end{align}
for almost all $g\in G$, where we identified $g=F_{Sp(N)}(\tilde{\phi}_1,\ldots,\xi)$ such that the parameters $\tilde{\phi}_1,\ldots,\xi$ are well-defined, $f^{SU(N)}_{ij}$ and $h^{SU(N)}_{ij}$ are finite-type functions on $SU(N)$, $m_{ij}, n_{ij}\in\Z$ and 
\begin{align*}
	g_{ij}(y_1,\ldots,y_N):=\sin^{p^1_{ij}}(y_1)\cos^{q^1_{ij}}(y_1)\ldots\sin^{p^N_{ij}}(y_N)\cos^{q^N_{ij}}(y_N)
\end{align*}
where $p_{ij}^n\in\mathbb{N}_0$ for all $n=1,2,\ldots,N$ and $q_{ij}^m\in\{0,1\}$ for all $m=1,2,\ldots,N$. We sum over two separate indices to make sure we have all the possible terms and allow for different powers of each term.

 For more details regarding the finite-type functions on $SU(N)$, we refer to \cite{Zwart}. In a similar way as in \cite{Zwart}, we have the following proposition 

\begin{proposition}\label{prop:transform_Mathieu_to_abelian_integral_Sp(N)}
	Let $f$ be a finite-type function on $Sp(N)$ as in Equation (\ref{eq:Finite_type_Sp(N)}). Then for any $P\in\N$ we have
	\begin{align*}
		\int_{Sp(N)}f^P(g) dg = C'&\int_{[0,1]^{N(N-1)}}\int_{(S^*)^{N(N+1)}}\int_{0}^1\int_0^{\xi_N}\ldots\int_0^{\xi_2}\hat{f}(x_1,\ldots,x_{N(N-1)},z_1,\ldots,z_{N(N+1)},\xi_1,\ldots,\xi_N)^P\\
		&\tilde{J}_{Sp(N)}(x_1,\ldots,x_{N(N-1)},\xi_1,\ldots,\xi_N)d\xi_1\ldots,d\xi_N \frac{dz_1}{z_1}\ldots \frac{dz_{N(N+1)}}{z_{N(N+1)}}dx_1\ldots dx_{N(N-1)}
	\end{align*}
	where
	\begin{align*}
		\hat{f}(x_1,\ldots,\xi_N):=\sum_{i,j}&c_{ij}\overline{f}_{ij}^{SU(N)}(x_1,\ldots,z_{\frac{N(N+1)}{2}-1})z_{\frac{N(N+1)}{2}}^{\frac{m_{ij}}{2}}\tilde{g}_{ij}(\xi_1,\ldots,\xi_N)\cdot\\
		&\widetilde{h_{ij}^{SU(N)}}(x_{\frac{N(N-1)}{2}+1},\ldots,z_{N(N+1)-1})z_{N(N+1)}^{n_{ij}},
	\end{align*}
	and 
	\begin{align*}
		\tilde{J}_{Sp(N)}(x_1,\ldots,x_{N(N-1)},\xi_1,\ldots,\xi_N):=\,&J_{SU(N)}(x_1,\ldots x_\frac{N(N-1)}{2})\left(\prod_{j=1}^N\xi_j\right)\cdot\\
		&\left(\prod_{j>k}\left(\xi_j^2(1-\xi_k^2)-(1-\xi_j^2)\xi_k\right)\right)\cdot\\
		&J_{SU(N)}(x_{\frac{N(N-1)}{2}+1},\ldots,x_{N(N-1)}).
	\end{align*}
	Here $C'>0$ is some constant and $\widetilde{h^{SU(N)}_{ij}}$ and $J_{SU(N)}$ are as in \cite[Lemma 2.7]{Zwart}. The function $J_{Sp(N)}$ is defined as in Theorem \ref{lemma:KAK_decomp_Sp(N)}, the set $S^*:=S^1\setminus\{1\}$ and the function $\overline{f}_{ij}^{SU(N)}$ is defined as $$\overline{f}^{SU(N)}_{ij}(x_1,\ldots, z_{\frac{N(N+1)}{2}-1})=\widetilde{f^{SU(N)}_{ij}}(x_1,\ldots,x_{\frac{N(N-1)}{2}},z_{1},\ldots,z_{\frac{N(N-1)}{2}},z_{\frac{N(N-1)}{2}+1}^{1/2},\ldots,z_{\frac{N(N+1)}{2}-1}^{1/2}),$$ and $$\tilde{g}_{ij}(\xi_1,\ldots,\xi_N)=\xi_1^{p^1_{ij}}(1-\xi_1^2)^{\frac{q_{ij}^1}{2}}\ldots \xi_N^{p_{ij}^N}(1-\xi_N^2)^{\frac{q_{ij}^N}{2}}.$$
\end{proposition}
\begin{proof}
	We note that
	\begin{align*}
		\int_G f(g)^P dg =& \int_G \left[\sum_{i,j} c_{ij}f^{SU(N)}_{ij}(\tilde{\phi_1},\ldots,\tilde{\omega}_{N-1})e^{im_{ij}\tilde{\xi}}\,g_{ij}(y_1,\ldots,y_N)\,h^{SU(N)}_{ij}(\phi_1\ldots,\omega_{N-1})e^{in_{ij}\xi}\right]^P dg\\
		=&\sum_{\sum\beta_{ij}=P}\int_G\prod_{i,j}c_{ij}^{\beta_{ij}}(f^{SU(N)}_{ij}(\tilde{\phi_1},\ldots,\tilde{\omega}_{N-1}))^{\beta_{ij}}e^{im_{ij}\beta_{ij}\tilde{\xi}}\,(g_{ij}(y_1,\ldots,y_N))^{\beta_{ij}}\cdot\\
		&\qquad\qquad(h^{SU(N)}_{ij}(\phi_1\ldots,\omega_{N-1}))^{\beta_{ij}}e^{in_{ij}\beta_{ij}\xi}dg\\
		=&\sum_{\sum\beta_{ij}=P}\left(\prod_{i,j}c_{ij}^{\beta_{ij}}\right)\left[C\int_{K/M}\prod_{i,j}(f^{SU(N)}_{ij}(\tilde{\phi}_1,\ldots,\tilde{\omega}_{N-1}))^{\beta_{ij}}e^{im_{ij}\beta_{ij}\tilde{\xi}}dg_K\right]\cdot\,\\
		&\left[\int_{\mathcal{A}}\prod_{i,j}(g_{ij}(y_1,\ldots,y_N))^{\beta_{ij}}J_{Sp(N)}(y_1,\ldots,y_N)dy_1\ldots dy_N\right]\cdot\\
		&\left[\int_{K}\prod_{i,j}(h^{SU(N)}_{ij}(\phi_1\ldots,\omega_{N-1}))^{\beta_{ij}}e^{in_{ij}\beta_{ij}\xi}dg_{K}\right].
	\end{align*}
	To complete the lemma, we consider each integral separately. Using \cite[Lemma 2.7]{Zwart}, we immediately can claim that the last integral equals to	
	\begin{align*}
		\int_{K}\prod_{i,j}&(h^{SU(N)}_{ij}(\phi_1\ldots,\omega_{N-1}))^{\beta_{ij}}e^{in_{ij}\beta_{ij}\xi}dg_{K} = \frac{1}{2(N-1)i^{\frac{N(N+1)}{2}}}\int_{[0,1]^{\frac{N(N-1)}{2}}}\int_{(S^*)^{\frac{N(N+1)}{2}}}\\
		&\prod_{i,j}(\widetilde{h_{ij}^{SU(N)}}(x_1,\ldots,z_{\frac{N(N+1)}{2}-1}))^{\beta_{ij}}z_{\frac{N(N+1)}{2}}^{n_{ij}\beta_{ij}}\cdot J_{SU(N)}(x_1,\ldots,x_{\frac{N(N-1)}{2}})\frac{dz_1}{z_1}\ldots \frac{dz_{\frac{N(N+1)}{2}}}{z_{\frac{N(N+1)}{2}}}dx_1\ldots dx_{\frac{N(N-1)}{2}},
	\end{align*}
	where $\widetilde{h_{ij}^{SU(N)}}$ is defined as in \cite[Lemma 2.7]{Zwart}, where we recall $S^*= S^1\setminus\{1\}$. The same arguments can be applied to the first integral, only note that the $\tilde{\omega}_j$ go over the range $[0,\frac{\pi}{j}]$ instead of $[0,\frac{2\pi}{j}]$ which gives a factor of $1/2$ per parameter. Therefore, we get that the first integral can be written as 
	\begin{align*}
	\int_{K}\prod_{i,j}&(f^{SU(N)}_{ij}(\tilde{\phi}_1\ldots,\tilde{\xi}))^{\beta_{ij}}e^{im_{ij}\beta_{ij}\tilde{\omega}_{N-1}}dg_{K} = \frac{1}{2^{N}(N-1)i^{\frac{N(N+1)}{2}}}\\
	&\int_{[0,1]^{\frac{N(N-1)}{2}}}\int_{(S^*)^{\frac{N(N+1)}{2}}}\prod_{i,j}(\overline{f}_{ij}^{SU(N)}(x_1,\ldots,z_{\frac{N(N+1)}{2}-1}))^{\beta_{ij}}z_{\frac{N(N+1)}{2}}^{\frac{m_{ij}}{2}\beta_{ij}}\cdot J_{SU(N)}(x_1,\ldots,x_{\frac{N(N-1)}{2}})\cdot\\
	& \qquad\qquad\qquad\qquad\qquad\qquad\frac{dz_1}{z_1}\ldots\frac{dz_{\frac{N(N+1)}{2}}}{z_{\frac{N(N+1)}{2}}}dx_1\ldots dx_{\frac{N(N-1)}{2}},
	\end{align*}
	where $$\overline{f}^{SU(N)}_{ij}(x_1,\ldots, z_{\frac{N(N+1)}{2}-1})=\widetilde{f^{SU(N)}_{ij}}(x_1,\ldots,x_{\frac{N(N-1)}{2}},z_{1},\ldots,z_{\frac{N(N-1)}{2}},z_{\frac{N(N-1)}{2}+1}^{1/2},\ldots,z_{\frac{N(N+1)}{2}-1}^{1/2}).$$
	Finally note that 
	\begin{align*}
		\int_{\mathcal{A}}\prod_{i,j}(g_{ij}(y_1,\ldots,y_N))^{\beta_{ij}}&J_{Sp(N)}(y_1,\ldots,y_N)dy_1\ldots dy_N\\
		= \int_0^{\pi/2}\int_0^{y_{N-1}}\ldots \int_0^{y_2}&\sin^{\sum_{ij}\beta_{ij}p_{ij}^1}(y_1)\cos^{\sum_{ij}\beta_{ij}q_{ij}^1}(y_1)\ldots\sin^{\sum_{ij}\beta_{ij}p_{ij}^N}(y_N)\cos^{\sum_{ij}\beta_{ij}q_{ij}^N}(y_N)\cdot\\
		&J_{Sp(N)}(y_1,\ldots,y_N)dy_1\ldots dy_N
	\end{align*}
	Note that
	\begin{align*}
		J(y_1,\ldots,y_N)&=\left[\prod_{j=1}^N\sin(2y_j)\right]\cdot\left[\prod_{j>k}\sin(y_j-y_k)\sin(y_j+y_k)\right]\\
		&=\left[\prod_{j=1}^N2\sin(y_j)\cos(y_j)\right]\cdot\left[\prod_{j>k}\left(\sin^2(y_j)\cos^2(y_k)-\cos^2(y_j)\sin^2(y_k)\right)\right].
	\end{align*}
	So if we subsitute $\xi_j=\sin(y_j)$, then we see that 
	\begin{align*}
		\int_{\mathcal{A}}&\prod_{i,j}(g_{ij}(y_1,\ldots,y_N))^{\beta_{ij}}J_{Sp(N)}(y_1,\ldots,y_N)dy_1\ldots dy_N\\
		=& \int_0^{1}\int_0^{\xi_{N-1}}\ldots \int_0^{\xi_2}\xi_1^{\sum_{ij}\beta_{ij}p_{ij}^1}(1-\xi_1^2)^{\frac{\sum_{ij}\beta_{ij}q_{ij}^1}{2}}\ldots \xi_N^{\sum_{ij}\beta_{ij}p_{ij}^N}(1-\xi_N^2)^{\frac{\sum_{ij}\beta_{ij}q_{ij}^N}{2}}\cdot\\
		& \left[\prod_{j=1}^N2\xi_j\right]\cdot\left[\prod_{j>k}\left(\xi_j^2(1-\xi_k^2)-(1-\xi_j^2)\xi_k\right)\right]d\xi_1\ldots d\xi_N
	\end{align*}
	Putting everything together gives:
	\begin{align*}
		\int_{Sp(N)} f^P(g)dg =&\sum_{\sum\beta_{ij}=P}\left(\prod_{i,j}c_{ij}^{\beta_{ij}}\right)\left[C\frac{1}{2^{N}(N-1)i^{\frac{N(N+1)}{2}}}\int_{[0,1]^{\frac{N(N-1)}{2}}}\int_{(S^*)^{\frac{N(N+1)}{2}}}\right.\\
		&\prod_{i,j}(\overline{f}_{ij}^{SU(N)}(x_1,\ldots,z_{\frac{N(N+1)}{2}-1}))^{\beta_{ij}}z_{\frac{N(N+1)}{2}}^{\frac{m_{ij}}{2}\beta_{ij}}\cdot J_{SU(N)}(x_1,\ldots,x_{\frac{N(N-1)}{2}})\cdot\\
		&\left.\frac{dz_1}{z_1}\ldots \frac{dz_{\frac{N(N+1)}{2}}}{z_{\frac{N(N+1)}{2}}}dx_1\ldots dx_{\frac{N(N-1)}{2}}\right]\cdot\\
		&\left[\int_0^{1}\int_0^{\xi_{N}}\ldots \int_0^{\xi_2}\xi_1^{\sum_{ij}\beta_{ij}p_{ij}^1}(1-\xi_1^2)^{\frac{\sum_{ij}\beta_{ij}q_{ij}^1}{2}}\ldots \xi_N^{\sum_{ij}\beta_{ij}p_{ij}^N}(1-\xi_N^2)^{\frac{\sum_{ij}\beta_{ij}q_{ij}^N}{2}}\right.\cdot\\
		&\left. \left(\prod_{j=1}^N2\xi_j\right)\cdot\left(\prod_{j>k}\left(\xi_j^2(1-\xi_k^2)-(1-\xi_j^2)\xi_k\right)\right)d\xi_1\ldots d\xi_N\right]\cdot\\
		&\left[\frac{1}{2(N-1)i^{\frac{N(N+1)}{2}}}\int_{[0,1]^{\frac{N(N-1)}{2}}}\int_{(S^*)^{\frac{N(N+1)}{2}}}\prod_{i,j}(\tilde{h}_{ij}^{SU(N)}(x_1,\ldots,z_{\frac{N(N+1)}{2}-1}))^{\beta_{ij}}z_{\frac{N(N+1)}{2}}^{n_{ij}\beta_{ij}}\cdot\right.\\
		&\left.J_{SU(N)}(x_1,\ldots,x_{\frac{N(N-1)}{2}})\frac{dz_1}{z_1}\ldots \frac{dz_{\frac{N(N+1)}{2}}}{z_{\frac{N(N+1)}{2}}}dx_1\ldots dx_{\frac{N(N-1)}{2}}\right]\\
		=&\frac{C}{2(N-1)^2 i^{N(N+1)}}  \int_{[0,1]^{N(N-1)}}\int_{(S^*)^{N(N+1)}}\int_{0}^1\int_0^{\xi_N}\ldots\int_0^{\xi_2}\\
		&\left[\sum_{i,j}c_{ij}\overline{f}_{ij}^{SU(N)}(x_1,\ldots,z_{\frac{N(N+1)}{2}-1})z_{\frac{N(N+1)}{2}}^{\frac{m_{ij}}{2}}\xi_1^{k_{ij}}(1-\xi_1^2)^{\frac{l_{ij}^1}{2}}\ldots \xi_N^{k_{ij}^N}(1-\xi_N^2)^{\frac{l_{ij}}{2}}\right.\cdot\\
		&\left.\tilde{h}_{ij}^{SU(N)}(x_{\frac{N(N-1)}{2}+1},\ldots,z_{N(N+1)-1})z_{N(N+1)}^{n_{ij}}\right]^PJ_{SU(N)}(x_1,\ldots x_\frac{N(N-1)}{2})\cdot\\
		&\left(\prod_{j=1}^N\xi_j\right)\cdot\left(\prod_{j>k}\left(\xi_j^2(1-\xi_k^2)-(1-\xi_j^2)\xi_k\right)\right) J_{SU(N)}(x_{\frac{N(N-1)}{2}+1},\ldots,x_{N(N-1)})\\
		&d\xi_1\ldots d\xi_N \frac{dz_1}{z_1}\ldots\frac{dz_{N(N+1)}}{z_{N(N-1)}}dx_1\ldots dx_{N(N-1)}.
	\end{align*}
	Defining 
	\begin{align*}
		\tilde{J}_{Sp(N)}(x_1,\ldots,&x_{N(N-1)},\xi_1,\ldots,\xi_N):=\,J_{SU(N)}(x_1,\ldots x_\frac{N(N-1)}{2})\left(\prod_{j=1}^N\xi_j\right)\cdot\\
		&\left(\prod_{j>k}\left(\xi_j^2(1-\xi_k^2)-(1-\xi_j^2)\xi_k\right)\right)J_{SU(N)}(x_{\frac{N(N-1)}{2}+1},\ldots,x_{N(N-1)}).
	\end{align*}
	and defining $C'=\frac{C}{(-1)^{\frac{N(N+1)}{2}}\cdot2(N-1)^2}$ gives the desired result.
\end{proof}


We note that Proposition \ref{prop:transform_Mathieu_to_abelian_integral_Sp(N)} gives a sum of products of (possible roots of) polynomials, in a similar way as in \cite{Zwart}. Let us give a name for these kind of polynomial

\begin{definition}\label{def:1/N_admissible}
	Let $k,l,N\in\N$ and $f: [0,1]^{k}\times (S^*)^{l}\rightarrow\C$. We say $f$ is a $\frac{1}{N}$-\emph{admissible function} if $f$ can be written as $$f(x_1,\ldots,x_k,z_1,\ldots,z_{l})=\sum_{\vec{m}}c_{\vec{m}}(x)z^{\vec{m}},$$ where $\vec{m}=(m_1,\ldots,m_l)$ is a multi-index where $m_i\in \bigcup_{j=1}^N\frac{1}{j}\mathbb{Z}$, and $c_{\vec{m}}(x)\in \C[x_1,(1-x_1^2)^{1/2},\ldots,x_k,(1-x_k^2)^{1/2}]$ is a complex polynomial in $x_i$ and $\sqrt{1-x_i^2}$. We define the collection of $\vec{m}$ for which $c_{\vec{m}}\neq 0$ \emph{the spectrum of} $f$, and it will be denoted by $\mathrm{Sp}(f)$.
\end{definition} 
\begin{remark}
	Note that in \cite{Zwart}, we called these functions $SU(N)$-admissible functions.
\end{remark}
It is clear from Proposition \ref{prop:transform_Mathieu_to_abelian_integral_Sp(N)} that $\hat{f}$ is a $\frac{1}{N}$-admissible function. Motivated by \cite{Mueger,Zwart}, we make the following conjecture:

\begin{conjecture}\label{con:xz-conjecture_Sp(N)}
	Let $f:[0,1]^{N^2}\times (S^*)^{N(N+1)}\rightarrow\C$ be a $\frac{1}{N}$-admissible function. If $$\int_{[0,1]^{N(N-1)}}\int_{(S^*)^{N(N+1)}}\int_0^1\int_0^{\xi_N}\ldots\int_0^{\xi_2}f^P \tilde{J}_{Sp(N)} \,\,d\xi_1\ldots d\xi_N\frac{dz_1}{z_1}\ldots \frac{dz_{N(N+1)}}{z_{N(N+1)}}dx_1\ldots dx_{N(N-1)} = 0$$ for all $P\in \N$, then $\vec{0}$ does not lie in the convex hull of $\mathrm{Sp}(f)$.
\end{conjecture}

\begin{theorem}
	Assume Conjecture \ref{con:xz-conjecture_Sp(N)} is true. Then the Mathieu conjecture is true for $Sp(N)$.
\end{theorem}
\begin{proof}
	Let $f,h$ be finite-type functions such that $\int_G f^P(g)dg=0$ for all $P\in\N.$ We are interested in $\int_G f^P(g)h(g)dg.$ Since $f,h$ are finite-type functions, by the previous discussion we can conclude $f$ and $h$ have the form of Equation (\ref{eq:Finite_type_Sp(N)}). By linearity of the integral we can assume $h$ to be of the form $$h(g)=c\cdot  h_1^{SU(N)}(\tilde{\phi}_1,\ldots,\tilde{\omega}_{N-1})e^{im\tilde{\xi}}h_2(y_1,\ldots,y_n)h_3^{SU(N)}(\phi_1,\ldots,\omega_{N-1})e^{in\xi}$$ where $c\in\C$, the functions $h_1, h_3$ are finite-type functions on $SU(N)$, $m,n\in\Z$ and $$h_2(y_1,\ldots,y_N)=\cos^{P_1}(y_1)\sin^{Q_1}(y_1)\ldots\cos^{P_N}(y_N)\sin^{Q_N}(y_N)$$ where $K_j\in \N_0$ and $L_j\in\{0,1\}$.
	Then the integral is of the form
	\begin{align}\label{eq:Mathieu_proof_Sp(N)_integral_that_is_nonzero}
		\begin{split}\int_G f^P(g)h(g)dg=&\sum_{\sum_{i,j}\beta_{ij}=P}\left(c\cdot\left(\prod_{i,j}c_{ij}^{\beta_{ij}}\right)\right.\\
		&\left[\int_{K/M}(f_{ij}^{SU(N)}(\tilde{\phi}_1,\ldots,\tilde{\omega}_{N-1})^{\beta_{ij}}e^{im_{ij}\beta_{ij}\tilde{\xi}}h_1^{SU(N)}(\tilde{\phi}_1,\ldots,\tilde{\omega}_{N-1})e^{im\tilde{\xi}}dg_{K/M}\right]\cdot\\
		&\left[\int_\mathcal{A}\prod_{i,j}(g_{ij}(y_1,\ldots,y_N))^{\beta_{ij}}h_2(y_1,\ldots,y_N)J_{Sp(N)}(y_1,\ldots,y_N)dy_1\ldots dy_N\right]\cdot\\
		&\left.\left[\int_{K}\prod_{i,j}(h_{ij}^{SU(N)}(\phi_1,\ldots,\omega_{N-1}))^{\beta_{ij}}e^{in_{ij}\xi}h_3^{SU(N)}(\phi_1,\ldots,\omega_{N-1})e^{in\xi}dg_{K}\right]\right).\end{split}
	\end{align}
	Assume to the contrary that $\int_G f(g)^Ph(g)dg\neq 0$ for infinitely many $P$. That means that there exists a set $\{\beta_{ij}\}_{i,j}$ such that the integrals in Equation (\ref{eq:Mathieu_proof_Sp(N)_integral_that_is_nonzero}) are non-zero. In particular, it means that 
	\begin{align}
		\int_{K/M}\prod_{i,j}(f_{ij}^{SU(N)}(\tilde{\phi}_1,\ldots,\tilde{\omega}_{N-1}))^{\beta_{ij}}e^{im_{ij}\beta_{ij}\tilde{\xi}}h_1^{SU(N)}(\tilde{\phi}_1,\ldots,\tilde{\omega}_{N-1})e^{im\tilde{\xi}}dg_{K/M}\neq 0\label{eq:Sp(N)_integral_over_U(N)_nr1}\\
		\int_\mathcal{A}\prod_{i,j}(g_{ij}(y_1,\ldots,y_N))^{\beta_{ij}}h_2(y_1,\ldots,y_N)\,J(y_1,\ldots,y_N)dy_1\ldots dy_N\neq 0 \label{eq:Sp(N)_integral_over_U(N)_nr2}\\
		\int_{K}\prod_{i,j}(h_{ij}^{SU(N)}(\phi_1,\ldots,\omega_{N-1}))^{\beta_{ij}}e^{in_{ij}\xi}h_3^{SU(N)}(\phi_1,\ldots,\omega_{N-1})e^{in\xi}dg_{K}\neq 0 \label{eq:Sp(N)_integral_over_U(N)_nr3}
	\end{align}
	Our goal now is to make extensive use of the Haar measure properties to conclude that $\vec{0}\in\mathrm{Sp}(f)$. Note that any element $g\in Sp(N)$ can be written as $$g=x\exp(H)k$$ as in Theorem \ref{thm:KAK_decomp}, where $x\in K/M$ is a representative of $xM$, $e^H\in A$ and $k\in K$, and the Haar measure can be written as $$dg_{Sp(N)}=dg_{K/M}\, J(\exp(H))dH\,dg_K.$$ Since $dg_K$ is just the Haar measure on $K$, for any smooth function $f:G\rightarrow\C$ we have 
	\begin{align*}
		\int_G f(g)dg &= \int_G f(x\exp(H)k)  |J(\exp(H))|\,dg_{K/M}\,dH\,dg_K\\
		&=\int_G f(x\exp(H)kh) |J(\exp(H))|\,dg_{K/M}\,dH\,dg_K\\
		&=\int_G f(x\exp(H)hk) |J(\exp(H))|\,dg_{K/M}\,dH\,dg_K
	\end{align*}
	for any $h\in K$. Choosing $k\in K$ we see that $kh, hk\in K$, and therefore by previous argument we see that only Equation (\ref{eq:Sp(N)_integral_over_U(N)_nr3}) would change. Also note that $dg_{K/M}$ is just the Haar measure on $K$ as well up to a factor. In the same way we thus see that 
	\begin{align*}
		\int_G f(g)dg &=\int_G f(xh\exp(H)k) |J(\exp(H))| \,dg_{K/M}\,dH\,dg_K\\
		&=\int_G f(hx\exp(H)k)|J(\exp(H))| \,dg_{K/M}\,dH\,dg_K
	\end{align*}
	where we used the left- and right-invariance of $dg_{K/M}$. Because of the diffeomorphism up to a measure zero in Theorem \ref{thm:KAK_decomp}, we have that only Equation (\ref{eq:Sp(N)_integral_over_U(N)_nr1}) would change. Being able to change Equation (\ref{eq:Sp(N)_integral_over_U(N)_nr1}) and Equation (\ref{eq:Sp(N)_integral_over_U(N)_nr3}) independently, we can conclude that Equations (\ref{eq:Sp(N)_integral_over_U(N)_nr1}) - (\ref{eq:Sp(N)_integral_over_U(N)_nr3}) are all independent equations.
	 
	But note that the Equation (\ref{eq:Sp(N)_integral_over_U(N)_nr1}) and Equation (\ref{eq:Sp(N)_integral_over_U(N)_nr3}) are statements of finite-type functions on $U(N)$. Expanding the finite-type functions on $SU(N)$, in the same line as in \cite{Zwart}, gives 
	\begin{align*}
		f_{ij}^{SU(N)}(\phi_1,\ldots,\omega_{N-1}):=&c_{ij} e^{ik^1_{ij}\phi_1}\sin^{a_{ij}^1}(\psi_1)\cos^{b_{ij}^1}(\psi_1)\ldots e^{ik^{N-1}_{ij}\phi_{N-1}}\sin^{a_{ij}^{N-1}}(\psi_{N-1})\cdot\\
		&\cos^{b_{ij}^{N-1}}(\psi_{N-1})f^{SU(N-1)}_{ij}(\phi_N,\ldots,\omega_{N-2})e^{il_{ij}^{N-1}\omega_{N-1}}
	\end{align*} and 
	\begin{align*}
		h_1^{SU(N)}(\phi_1,\ldots,\omega_{N-1})=&c_1 e^{iK_1\phi_1}\sin^{A_1}(\psi_1)\cos^{B_1}(\psi_1)\ldots e^{iK_{N-1}\phi_{N-1}}\sin^{A_{N-1}}(\psi_{N-1})\cdot\\
		&\cos^{B_{N-1}}(\psi_{N-1})(h_{SU(N-1)})(\phi_N,\ldots,\omega_N)e^{iL_{N-1}\omega_{N-1}}
	\end{align*} where $f_{ij}^{SU(N-1)}$ and $h_{SU(N-1)}$ are defined inductively, and the indices satisfy
	\begin{align*}
		k_{ij}^1,\ldots,k_{ij}^{N-1},k_{ij}^{N},\ldots,k_{ij}^{N(N-1)/2},l_{ij}^1,\ldots,l_{ij}^{N-1},K_1,\ldots,K_{N(N-1)/2},L_1,\ldots,L_{N-1}\in\Z,\\
		a_{ij}^1,\ldots a_{ij}^{N(N-1)/2},A_1,\ldots,A_{N(N-1)}\in\N_0,\qquad b_{ij}^1,\ldots b_{ij}^{N(N-1)/2},B_1,\ldots,B_{N(N-1)/2}\in\{0,1\}.
	\end{align*}
	Applying the same steps as in the proof of \cite[Thm. 2.11]{Zwart}, we find that $$\sum_{i,j}\beta_{ij}k_{ij}^m+K_m=0,\qquad\sum_{i,j}\beta_{ij}l_{ij}^n+L_n=0$$ for all $m=1,\ldots,\frac{N(N-1)}{2}$ and $n=1,\ldots,N$. For more details we refer to \cite{Zwart}. Using the same arguments for $\tilde{f}_{ij}^{SU(N)}$, assigning a tilde on the indices of $\tilde{f}_{ij}^{SU(N)}$, we find that
	$$\sum_{i,j}\beta_{ij}\tilde{k}_{ij}^m+\tilde{K}_m=0,\qquad\sum_{i,j}\beta_{ij}\tilde{l}_{ij}^n+\tilde{L}_n=0$$ for all $m=1,\ldots,\frac{N(N-1)}{2}$ and $n=1,\ldots,N$.
	Finally, we note that $$\int f(xe^Hk)^P |J(\exp(H))|dk_M\,dH\,dk = \int f\left(xe^Hk\exp\begin{pmatrix}
	t\lambda_0 & 0\\
	0&0
	\end{pmatrix}\right)|J(\exp(H))|dk_M\,dH\,dk,$$ which must be equivalent to the map $\tilde{\xi}\mapsto \tilde{\xi}+t$. This mapping can only be invariant for all $t\in\R$ if $$\sum_{ij}n_{ij}\beta_{ij}+n=0.$$ In the same way we get 
	$$\int f(xe^Hk)^P |J(\exp(H))|dk_M\,dH\,dk = \int f\left(x\exp\begin{pmatrix}
	t\lambda_0 & 0\\
	0&0
	\end{pmatrix}e^Hk\right)|J(\exp(H))|dk_M\,dH\,dk,$$ which gives in the same way $$\sum_{i,j}m_{ij}\beta_{ij}+m=0.$$
	In other words, adding it all together, we get that $$\sum_{ij}(\sum_mk_{ij}^m + \tilde{k}_{ij}^m+n_{ij}+m_{ij})\beta_{ij}=-\left(\sum_{m}(K_m+\tilde{K}_m)+n+m\right).$$ Multiplying both sides with $\frac{1}{P}$ and letting $P$ go to zero, gives a sequence that converges to 0. Since the convex hull is a closed set, $0\in \mathrm{Conv}(\mathrm{Sp}(f))$. But this is a contradiction with Conjecture \ref{con:xz-conjecture_Sp(N)}. So $\int_{Sp(N)}f(g)^Ph(g)dg\neq0$ for only finitely many $P$, so $\int_{Sp(N)}f(g)^Ph(g)dg=0$ for $P$ large enough.
\end{proof}

\section{The case of $G_2$}\label{sec:G2}
\subsection{The Euler angle decomposition of $G_2$}
Next, we shift our focus to $G_2$. The structure of this section will be the same as before: we first find an Euler angles decomposition and then translate the Mathieu conjecture to a conjecture on $[0,1]^m\times (S^1)^n$ for some $n,m\in\N$. However, $G_2$ is not by definition defined as a closed subgroup of $GL(n,\C)$ for some $n$. Although that is not a necessity to apply Theorem \ref{thm:KAK_decomp}, we will do so to make much of the calculations actually computable. For that, we make extensive use of \cite{Cacciatori2} which gives an explicit embedding of $G_2$ into $GL(7,\R)$ and a way to do the Euler angles decomposition. 

The group $G_2$ can be seen as the set of automorphisms on the octonions \cite{OctonionsBaez}. Let $\mathbb{O}$ be the set of octonions. We can see $\mathbb{O}$ as a 8-dimensional vectorspace over $\R$, the spanning vectors being $e_0,\ldots,e_7$ where $e_0$ corresponds to the real unit, and $e_j$ to imaginary units for all $j=1,\ldots,7$. There exists a natural multiplication on $\mathbb{O}$ which we will not need for our purpose. For more information on the octonions we refer to \cite{OctonionsBaez}. However, it is known \cite{OctonionsBaez} that $G_2$ can be seen as the automorphism group of $\mathbb{O}$. Since $A$ is invertible and is linear, it must leave $e_0$ fixed and can only permute $e_1,\ldots,e_7$. Therefore one can consider $G_2$ as a closed subset of $GL(7,\R)$.

Since the Lie algebra of automorphisms on an algebra $A$ is the set of derivations on $A$, the Lie algebra $\mathfrak{g}_2\subset \mathfrak{gl}(7,\R)$ should therefore be the set of derivations on $\mathbb{O}$. One can compute these, and the generators are given in Appendix \ref{appendix:spanning_set_G_2}.

We note that $G_2$ is a simply connected compact simple Lie group, hence we can apply Theorem \ref{thm:KAK_decomp}. In the same way as in $Sp(N)$, we first find a Cartan involution. We choose the analytic involutive automorphism given by $$\theta=\Ad\begin{pmatrix}
\mathbf{1}_3&\\
&-\mathbf{1}_4
\end{pmatrix}.$$ It is clear that $\theta\in \mathrm{Int}(\mathfrak{g}_2)$ and one can easily see that only matrices of the form $\begin{pmatrix}
A&0\\
0&B
\end{pmatrix}$, where $A\in \mathfrak{gl}(3,\R)$ and $B\in \mathfrak{gl}(4,\R)$, are invariant under $\theta$. The Lie algebra $\mathfrak{k}$ is then given by $$\mathfrak{k}:=\mathrm{span}_{\R}(\lambda_{1},\lambda_2,\lambda_3,\lambda_8,\lambda_9,\lambda_{10}) = \mathfrak{k}_1\oplus \mathfrak{k}_2$$ where $\mathfrak{k}_1:=\mathrm{span}_\R(\lambda_1,\lambda_2,\lambda_3)$ and $\mathfrak{k}_2:=\mathrm{span}_\R(\lambda_8,\lambda_9,\lambda_{10})$, both isomorphic to the simple Lie algebra $\mathfrak{su}(2)$ and $[\mathfrak{k}_1,\mathfrak{k}_2]=0$. Note that the matrices $\{\lambda_j\}_j$ are defined in Appendix \ref{appendix:spanning_set_G_2}. This defines our $K$ to be $$K:=\langle e^\mathfrak{k}\rangle.$$ Since $[\mathfrak{k}_1,\mathfrak{k}_2]=0$ one would expect $K$ to be isomorphic to $SU(2)\times SU(2)$. Let us define the latter set by $$\hat{K}:=\langle e^{\mathfrak{k}_1}\rangle\times \langle e^{\mathfrak{k}_2}\rangle \simeq SU(2)\times SU(2).$$ However, we note that $K\neq \hat{K}$ since the group homomorphism $$f:\hat{K}\rightarrow GL(7,\R),\qquad f(x,y)=xy$$ has $\ker{f}=\{(\mathbf{1}_7,\mathbf{1}_7),(\mathrm{diag}(\mathbf{1}_3,-\mathbf{1}_4),\mathrm{diag}(\mathbf{1}_3,-\mathbf{1}_4))\}$. So we see that $$K \simeq (SU(2)\times SU(2))/\mathbb{Z}_2\simeq SO(4)$$ where $\Z_2\simeq \ker{f}$ and $\hat{K}$ is the double cover of $K$.

So to apply Theorem \ref{thm:KAK_decomp}, let us first describe the Euler angle decomposition of $K$. The Euler angle decomposition on $\hat{K}$ can be described as follows

\begin{lemma}
	Let $\hat{K}$ be as before, then up to a measure zero set there is a diffeomorphism such that $$\hat{K} \simeq K_1A_1K'_1\times K_2A_2K'_2$$ as manifolds. Here 
	\begin{align*}
		K_1&:=\{e^{x\lambda_3}\,|\,x\in [0,\pi]\},\qquad\,\, K_2:=\{e^{x\lambda_{8}}\,|\,x\in [0,\pi]\}\\
		A_1&:=\left\{e^{x\lambda_2}\left|\,x\in \left[0,\frac{\pi}{2}\right.\right]\right\},\quad A_2:=\left\{e^{x\lambda_9}\left|\,x\in \left[0,\frac{\pi}{2}\right]\right\}\right.\\
		K'_1&:=\{e^{x\lambda_3}\,|\,x\in [0,2\pi]\},\qquad  K'_2:=\{e^{x\lambda_{8}}\,|\,x\in [0,2\pi]\}.
	\end{align*}
\end{lemma}
\begin{proof}
	The fact that $\C$ can be seen as $\R^2$, where multiplication with $i$ can seen as applying the matrix 
	\begin{align*}
	\begin{pmatrix}
		0&1\\
		-1&0
	\end{pmatrix}
	\end{align*}
	we see that $SU(2)$ can be readily been embedded into $SO(4)$, and we know that $SU(2)$ can be described by the Euler angle decomposition as in, for example, \cite{Dings-Koelink,Mueger,Zwart} by
	\begin{align*}
		SU(2)\simeq K_1AK_2
	\end{align*} 
	where $$K_1=\left\{\left.\exp\begin{pmatrix}ix&0\\ 0&-ix\end{pmatrix}\right|x\in[0,\pi)\right\},\qquad A=\left\{\left.\exp\begin{pmatrix}0&x\\ -x&0\end{pmatrix}\right|x\in[0,\pi)\right\}$$ and $$K_2=\left\{\left.\exp\begin{pmatrix}ix&0\\ 0&-ix\end{pmatrix}\right|x\in[0,2\pi)\right\}.$$ Applying this to every component gives the result.
\end{proof}
This way, we see that if $k\in \hat{K}$ then
\begin{align*}
	k = (e^{\phi_1 \lambda_3}e^{\psi_1\lambda_2}e^{\omega_1 \lambda_3},e^{\phi_2 \lambda_{8}}e^{\psi_2 \lambda_9}e^{\omega_2\lambda_{8}})
\end{align*}
where $\phi_1,\phi_2\in[0,\pi]$, $\psi_1,\psi_2\in[0,\pi/2]$ and $\omega_1,\omega_2\in [0,2\pi]$. Note that $e^{\pi \lambda_3}=\mathrm{diag}(\mathbf{1}_3,-\mathbf{1}_4)$, and since $e^{\mathfrak{k}_1}e^{\mathfrak{k}_2}=e^{\mathfrak{k}_2}e^{\mathfrak{k}_1}$, we can pull the mod $\mathbb{Z}_2$ operation into the first component, i.e. we get a unique description of $g\in K$ up to a measure zero set if we restrict the range of $\omega_1$ to $\omega_1\in[0,\pi]$. So we see that the following lemma holds

\begin{lemma}[Euler angle decomposition of $K$]\label{lemma:Parametrization_K}
	Let $K$ be as before. Define the mapping $F_K:[0,\pi]\times\left[0,\frac{\pi}{2}\right]\times [0,\pi]^2\times [0,\frac{\pi}{2}]\times [0,2\pi]\rightarrow K$ by
	\begin{align*}
		F_K(\phi_1,\psi_1,\omega_1,\phi_2,\psi_2,\omega_2):= e^{\phi_1 \lambda_3}e^{\psi_1\lambda_2}e^{\omega_1 \lambda_3}e^{\phi_2 \lambda_{8}}e^{\psi_2 \lambda_9}e^{\omega_2\lambda_{8}}.
	\end{align*}
	Then this mapping is a smooth surjective map, and is a diffeomorphism up to a measure zero set.
\end{lemma}
With the Euler angle decomposition of $K$ found, we note that the pair $(G_2,K)$ is the Riemannian symmetric pair associated with $(\mathfrak{g}_2,\theta)$. To get the rest of the data we need for Theorem \ref{thm:KAK_decomp}, we need to find a maximal abelian subalgebra $\mathfrak{a}\subseteq \mathfrak{p}$, where $$\mathfrak{p}:=\mathrm{span}_\R(\lambda_4,\lambda_5,\lambda_6,\lambda_7,\lambda_{11},\lambda_{12},\lambda_{13},\lambda_{14}).$$ We choose the maximal abelian subalgebra to be $$\mathfrak{a}:=\mathrm{span}_\R(\lambda_5,\lambda_{11}).$$
To apply Theorem \ref{thm:KAK_decomp}, we need to calculate $M$. By definition, $x\in M=Z_{K}(\mathfrak{a})$ if and only if $x\in K$ and $$\Ad(x)H=H$$ for all $H\in \mathfrak{a}$. A direct computation shows that, if we define
\begin{align*}
	\sigma:=\scalemath{0.8}{\begin{pmatrix}
			1&0&0&0&0&0&0\\
			0&-1&0&0&0&0&0\\
			0&0&-1&0&0&0&0\\
			0&0&0&1&0&0&0\\
			0&0&0&0&1&0&0\\
			0&0&0&0&0&-1&0\\
			0&0&0&0&0&0&-1
	\end{pmatrix}},\qquad \eta:=\scalemath{0.8}{\begin{pmatrix}
			-1&0&0&0&0&0&0\\
			0&0&1&0&0&0&0\\
			0&1&0&0&0&0&0\\
			0&0&0&-1&0&0&0\\
			0&0&0&0&1&0&0\\
			0&0&0&0&0&0&-1\\
			0&0&0&0&0&-1&0
	\end{pmatrix}},
\end{align*}
then $M=\{1,\sigma,\eta,\sigma\eta\}$.  With that, we get the following lemma

\begin{lemma}[Euler Angle decomposition of $G_2$]\label{lemma:Euler_angles_G2}
	Consider the compact simple Lie group $G_2$. Let $K=\langle e^{\mathfrak{k}}\rangle$ be as before. Then $$G_2\simeq (K/M)\,A\, K$$ diffeomorphic up to a measure zero set, where $A=\exp\mathcal{A}$ and $M=\{1,\sigma,\eta,\sigma\eta\}$ as before. In addition, define the mapping $$F_{G_2}:[0,\pi]\times \left[0,\frac{\pi}{4}\right]\times \left[0,\frac{\pi}{2}\right]\times \left[0,\pi\right]\times \left[0,\frac{\pi}{2}\right]\times [0,2\pi]\times\mathcal{A}\times [0,\pi]\times \left[0,\frac{\pi}{2}\right]\times[0,\pi]^2\times\left[0,\frac{\pi}{2}\right]\times[0,2\pi]\rightarrow G_2$$ by 
	\begin{align*}
		F_{G_2}&(\tilde{\phi}_1,\tilde{\psi}_1,\tilde{\omega}_1,\tilde{\phi}_2,\tilde{\psi}_2,\tilde{\omega}_2,y_1,y_2,\phi_1,\psi_1,\omega_1,\phi_2,\psi_2,\omega_2)=\\
		& e^{\tilde{\phi}_1 \lambda_3}e^{\tilde{\psi}_1\lambda_2}e^{\tilde{\omega}_1 \lambda_3}e^{\tilde{\phi}_2 \lambda_8}e^{\tilde{\psi}_2 \lambda_9}e^{\tilde{\omega}_2\lambda_8} e^{y_1\lambda_5}e^{y_2\lambda_{11}}e^{\phi_1 \lambda_3}e^{\psi_1\lambda_2}e^{\omega_1 \lambda_3}e^{\phi_2 \lambda_8}e^{\psi_2 \lambda_9}e^{\omega_2\lambda_8}
	\end{align*}
	where we repeat $\tilde{\phi}_1,\tilde{\phi}_2,\phi_1,\omega_1,\phi_2\in[0,\pi]$, $\tilde{\psi}_1\in[0,\frac{\pi}{4}]$, $\tilde{\omega}_1,\tilde{\psi}_2,\psi_1,\psi_2\in[0,\frac{\pi}{2}]$ and $\tilde{\omega}_2,\omega_2\in[0,2\pi]$ and $y_1,y_2\in \mathcal{A}$ for clarity.
	Here $$\mathcal{A}=\left\{(y_1,y_2)\in\R^2\left|0\leq y_i\leq \frac{\pi}{2}\text{ for } i=1,2\text{ and such that }y_2\leq \frac{1}{3}y_1 \right\}\right..$$
	Then $F_{G_2}$ is surjective and, up to a measure zero set, a diffeomorphism onto $G_2$.
	
	Finally, the Haar measure $dg$ decomposes as follows
	\begin{align*}
		dg = CJ_{G_2}(y_1, y_2) dg_{K/M} dy_1dy_2 dg_K 
	\end{align*}
	where $C>0$ is some constant, $dg_K$ is the Haar measure on $K$, $dg_{K/M}$ the unique left-invariant measure on $K/M$, both given by
	\begin{align}
		dg_K&=\cos(\psi_1)\sin(\psi_1)\cos(\psi_2)\sin(\psi_2)\,d\phi_1\,d\psi_1\,d\omega_1\, d\phi_2\, d\psi_2\,d\omega_2,\\
		dg_{K/M}&=  \cos(\tilde{\psi}_1)\sin(\tilde{\psi}_1)\cos(\tilde{\psi}_2)\sin(\tilde{\psi}_2)\,d\tilde{\phi}_1\,d\tilde{\psi}_1\,d\tilde{\omega}_1\,d\tilde{\phi}_2\, d\tilde{\psi}_2\, d\tilde{\omega}_2
	\end{align} and 
	\begin{align*}
		J_{G_2}(y_1,y_2):= \sin(y_1-3y_2)\sin(y_1-y_2)\sin(y_1+y_2)\sin(y_1+3y_2)\sin(2y_1)\sin(2y_2).
	\end{align*}
	Here the integration over the $y_j$ variables goes exactly over the area $\mathcal{A}$.
\end{lemma}
\begin{proof}
	The first part of the Theorem is proven by just applying Theorem \ref{thm:KAK_decomp} to $G_2$. To prove the claims for $F_{G_2}$, we need to find a parametrization on $K$, on $A$ and on $K/M$. The parametrization on $K$ is given by Lemma \ref{lemma:Parametrization_K}, and the parametrization of $A$ is immediate since $\mathfrak{a}$ is abelian. So we only need to find a parametrization of $K/M$. We recall that $M$ is discrete and finite, so we can describe $K/M$ as a subset of $K$. So let $g\in K$. Using Lemma \ref{lemma:Parametrization_K} to get the decomposition $g=e^{\phi_1 \lambda_3}e^{\psi_1\lambda_2}e^{\omega_1 \lambda_3}e^{\phi_2 \lambda_8}e^{\psi_2 \lambda_9}e^{\omega_2\lambda_8}$, we see immediately that 
	\begin{align*}
		g\sigma = e^{\phi_1 \lambda_3}e^{\psi_1\lambda_2}e^{(\omega_1+\frac{\pi}{2}) \lambda_3}e^{\phi_2 \lambda_8}e^{\psi_2 \lambda_9}e^{(\omega_2+\frac{\pi}{2})\lambda_8}.
	\end{align*}
	In the same way, after some calculations, we see that $$g\eta = e^{\phi_1 \lambda_3}e^{(\psi_1+\frac{\pi}{2})\lambda_2}e^{(-\frac{\pi}{4}-\omega_1) \lambda_3}e^{\phi_2 \lambda_8}e^{(\psi_2+\frac{\pi}{2}) \lambda_9}e^{(\frac{\pi}{4}-\omega_2)\lambda_8}.$$
	Note we need to decompose this more, for in our decomposition $\psi_j\in[0,\frac{\pi}{2}]$ and $g\eta$ is not described by this range. However, note that $$e^{\psi_1\lambda_2}=e^{\frac{\pi}{2}\lambda_3}e^{(\pi-\psi_1)\lambda_2}e^{\frac{\pi}{2}\lambda_3}$$ and the similarly for $e^{\psi_2\lambda_9}$. In other words, we see that
	$$g\eta = e^{(\phi_1+\frac{\pi}{2}) \lambda_3}e^{(\frac{\pi}{2}-\psi_1)\lambda_2}e^{(\frac{\pi}{4}-\omega_1) \lambda_3}e^{(\phi_2+\frac{\pi}{2}) \lambda_8}e^{(\frac{\pi}{2}-\psi_2) \lambda_9}e^{(\frac{3\pi}{4}-\omega_2)\lambda_8}.$$
	Therefore, the equivalence relation on $K$ to go to $K/M$ is given by 
	\begin{align*}
		g = g\sigma \qquad&\Leftrightarrow\qquad(\phi_1,\psi_1,\omega_1,\phi_2,\psi_2,\omega_2)\sim(\phi_1,\psi_1,\omega_1+\frac{\pi}{2},\phi_2,\psi_2,\omega_2+\frac{\pi}{2})\\
		g = g\eta \qquad&\Leftrightarrow\qquad(\phi_1,\psi_1,\omega_1,\phi_2,\psi_2,\omega_2)\sim\left(\phi_1+\frac{\pi}{2},\frac{\pi}{2}-\psi_1,\frac{\pi}{4}-\omega_1,\phi_2+\frac{\pi}{2},\frac{\pi}{2}-\psi_2,\frac{\pi}{4}-\omega_2\right)
	\end{align*}
	So if we reduce our interval of $\omega_1$ to $\omega_1\in\left[0,\frac{\pi}{2}\right]$ and the interval of $\psi_1$ to $\psi_1\in \left[0,\frac{\pi}{4}\right]$, we can describe $kM\in K/M$ explicitly by $$kM=e^{\phi_1 \lambda_3}e^{\psi_1\lambda_2}e^{\omega_1 \lambda_3}e^{\phi_2 \lambda_8}e^{\psi_2 \lambda_9}e^{\omega_2\lambda_8}$$ for some unique $\phi_1,\phi_2,\in[0,\pi]$, $\psi_1\in[0,\frac{\pi}{4}]$, $\omega_1,\psi_2\in[0,\frac{\pi}{2}]$ and $\omega_2\in[0,2\pi]$, up to a measure zero set.
	
	Now we are in the position to apply Theorem \ref{thm:KAK_decomp}, by applying both our Euler angle decompositions of $K$, as given in Lemma \ref{lemma:Parametrization_K}, and $K/M$ which was given above. This gives that any $g\in G_2$ can be written uniquely, up to a measure zero set, as
	\begin{align*}
		g&=k_1ak_2\\
		&= e^{\tilde{\phi}_1 \lambda_3}e^{\tilde{\psi}_1\lambda_2}e^{\tilde{\omega}_1 \lambda_3}e^{\tilde{\phi}_2 \lambda_8}e^{\tilde{\psi}_2 \lambda_9}e^{\tilde{\omega}_2\lambda_8} e^{y_1\lambda_5}e^{y\lambda_{11}}e^{\phi_1 \lambda_3}e^{\psi_1\lambda_2}e^{\omega_1 \lambda_3}e^{\phi_2 \lambda_8}e^{\psi_2 \lambda_9}e^{\omega_2\lambda_8}
	\end{align*}	
	where $\tilde{\phi}_1,\tilde{\phi}_2,\phi_1,\omega_1,\phi_2\in[0,\pi]$, $\tilde{\psi}_1\in[0,\frac{\pi}{4}]$, $\tilde{\omega}_1,\tilde{\psi}_2,\psi_1,\psi_2\in[0,\frac{\pi}{2}]$ and $\tilde{\omega}_2,\omega_2\in[0,2\pi]$ and $y_1,y_2\in\mathcal{A}$, where $\mathcal{A}$ is as in Theorem \ref{thm:KAK_decomp}.
	
	Note that both $K$ and $K/M$ can be seen as open subsets of $SU(2)\times SU(2)$, hence the Haar measure is given by the Haar measure on $SU(2)\times SU(2)$ restricted to $K$ and $K/M$ times some constant, i.e. there exist $C',C''>0$ such that the Haar measure on $K$ and the left-invariant measure on $K/M$ are given by $$C' dg_K=C''dg_{K/M}=dg_{SU(2)\times SU(2)}=\cos(\psi_1)\sin(\psi_1)\cos(\psi_2)\sin(\psi_2)\, d\phi_1\,d\phi_2d\psi_1\,d\psi_2\,d\omega_1\,d\omega_2$$ proving the Haar measure decomposition. To get an explicit form of $\mathcal{A}$, and the corresponding Haar measure on $\exp{\mathfrak{\mathcal{A}}}$, we recall the root decomposition of $\mathfrak{g}_2$ with respect to $\mathfrak{a}$. The details of the root decomposition are in Appendix \ref{appendix:spanning_set_G_2}. Recall $\mathfrak{a}=\mathrm{span}(\lambda_5,\lambda_{11})$. So we define the linear functionals $\alpha,\beta\in\mathfrak{a}^*$ as 
	\begin{align*}
		\alpha(\lambda_5)&:=0,\qquad \alpha(\lambda_{11}):=2\\
		\beta(\lambda_5)&:=1,\qquad \beta(\lambda_{11}):=-3.
	\end{align*}
	As will be shown in Appendix \ref{appendix:spanning_set_G_2}, the functionals $i\alpha,i\beta$ are the simple roots, and the set of positive roots are given by $$\Delta_\mathfrak{p}^+=\{i\alpha,i\beta,i(\alpha+\beta),i(2\alpha+\beta),i(3\alpha+\beta),i(3\alpha+2\beta)\}.$$
	It is then immediate from Equation (\ref{eq:Jacobian_abstract_Euler_decomp}) that for $H=y_1\lambda_5+y_2\lambda_{11}$ we get
	\begin{align*}
		J_{G_2}(\exp(H))=&\sin(2y_2)\sin(y_1-3y_2)\sin(y_1+2y_2-3y_2)\sin(y_1+4y_2-3y_2)\sin(y_1+6y_2-3y_2)\cdot\\
		&\sin(2y_1+6y_1-6y_1)\\
		=&\sin(2y_2)\sin(2y_1)\sin(y_1-3y_2)\sin(y_1-y_2)\sin(y_1+y_2)\sin(y_1+3y_2)\sin(2y_2).
	\end{align*}
	With that, we see that $$\mathcal{A}=\{(y_1,y_2)\in\R^2\,|\,0\leq2y_i\leq\pi \text{ for }i=1,2\text{ and such that }0\leq y_1 \pm y_2\leq\pi\text{ and }0\leq y_1 \pm 3y_2\leq\pi\}.$$ Note that this means that $0\leq y_i\leq \frac{\pi}{2}$, hence $y_1\pm y_2\leq \pi$ does not give any new information. In the same way we see that $0\leq y_1-3y_2$ gives that $y_2\leq \frac{1}{3}y_1$ which is also satisfied by $0\leq y_1 - y_2$. Finally if $y_1\leq \frac{\pi}{2}$ then $y_1+3y_2\leq \frac{\pi}{2}+\frac{\pi}{2}$. In other words, we can equally describe $\mathcal{A}$ by
	\begin{align}
		\mathcal{A}=\left\{(y_1,y_2)\in\R^2\left|0\leq y_i\leq \frac{\pi}{2}\text{ for } i=1,2\text{ and such that }y_2\leq \frac{1}{3}y_1 \right\}\right.
	\end{align}
	which concludes the proof.
\end{proof}

\subsection{The Mathieu Conjecture on $G_2$}\label{sec:Mathieu_G_2}
In this section, we will follow the same reasoning as in Section \ref{sec:Mathieu_G_2} to reduce the Mathieu conjecture to a more abelian version. We note that $G_2\subset SL(7,\C)$, hence the finite-type functions of $G_2$ are generated by the matrix entries and the inverse determinant. Since we have, up to a measure zero set, a description of the matrix entries using Lemma \ref{lemma:Euler_angles_G2}, we see that any finite type function $f$ is given by, using the diffeomorphism $g=F_{G_2}(\tilde{\phi}_1,\ldots,\omega_2)$,
\begin{align}
	\begin{split}\label{eq:finite_type_f_G2}
	f(g)=\sum_{j=1}^M\sum_{i=1}^Q& c_{ij}e^{ik_{ij}^1\tilde{\phi}_1}\sin^{l_{ij}^1}(\tilde{\psi}_1)\cos^{m_{ij}^1}(\tilde{\psi}_1) e^{ik_{ij}^2\tilde{\omega}_1}e^{ik_{ij}^3\tilde{\phi}_2}\sin^{l_{ij}^3}(\tilde{\psi}_2)\cos^{m_{ij}^3}(\tilde{\psi}_2)\cdot\\
	&e^{ik_{ij}^4\tilde{\omega}_2}\sin^{l_{ij}^4}(y_1)\cos^{m_{ij}^4}(y_1)\sin^{l_{ij}^5}(y_2)\cos^{m_{ij}^5}(y_2)e^{ik_{ij}^5\phi_1}\sin^{l_{ij}^6}(\psi_1)\cos^{m_{ij}^6}(\psi_1)\cdot\\
	 &e^{ik_{ij}^6\omega_1}e^{ik_{ij}^7\phi_2}\sin^{l_{ij}^{7}}(\psi_2)\cos^{m_{ij}^{7}}(\psi_2)e^{i k_{ij}^8\omega_2}.
	\end{split}
\end{align}
Here $k_{ij}^p\in \Z$, $l_{ij}^p\in \N_0$ and $m_{ij}^p\in\{0,1\}$ for all $i,j,p$. We reduced $m_{ij}^p$ to be in $\{0,1\}$ by using $\sin(x)^2+\cos(x)^2=1$. In a similar way as in Section \ref{sec:Mathieu_Sp(N)}, we can rewrite the condition $\int_{G_2}f^P\,dg$ in the following way

\begin{proposition}\label{prop:transform_Mathieu_to_abelian_integral_G_2}
	Let $f$ be a finite-type function on $G_2$ as in Equation (\ref{eq:finite_type_f_G2}). Then for any $P\in\N$ we have
	\begin{align*}
		\int_{G_2}f^P(g)dg = C'\int_{(S^*)^{8}}\int_{[0,1]^5}\int_{0}^{S(\xi_1)}\tilde{f}_{G_2}^P \tilde{J}_{G_2}(x_1,\ldots,x_4,\xi_1,\xi_2)d\xi_1d\xi_2 dx_1\ldots dx_4\frac{dz_1}{z_1}\ldots\frac{dz_8}{z_8}
	\end{align*}
	where $C'=\frac{C}{2^5\sqrt{2}}>0$ is some constant, 
	\begin{align*}
		\tilde{f}_{G_2}(z_1,\ldots,z_7,x_1,\ldots,x_4,\xi_1,\xi_2):=\sum_{j=1}^M\sum_{i=1}^Q& c_{ij}z_1^{\frac{k_{ij}^1}{2}}\left(\frac{x_1}{\sqrt{2}}\right)^{l_{ij}^1}\left(1-\frac{x_1^2}{2}\right)^{\frac{m_{ij}^1}{2}} z_2^{\frac{k_{ij}^2}{4}}z_3^{\frac{k_{ij}^3}{2}}x_2^{l_{ij}^3}(1-x_2^2)^{\frac{m_{ij}^3}{2}}\cdot\\
		&z_4^{k_{ij}^4}\xi_1^{l_{ij}^4}(1-\xi_1^2)^{\frac{m_{ij}^4}{2}}\xi_2^{l_{ij}^5}(1-\xi_2^2)^{\frac{m_{ij}^5}{2}}z_5^{\frac{k_{ij}^5}{2}}x_3^{l_{ij}^6}(1-x_3^2)^{\frac{m_{ij}^6}{2}}\cdot\\
		&z_6^{\frac{k_{ij}^6}{2}}z_7^{\frac{k_{ij}^{7}}{2}}x_4^{l_{ij}^{7}}(1-x_4^2)^{\frac{m_{ij}^{7}}{2}}z_8^{k_{ij}^8},
	\end{align*} and 
	\begin{align}
	\begin{split}
		\tilde{J}_{G_2}(x_1,\ldots,x_4,\xi_1,\xi_2):=& \xi_1\xi_2\bigg[\xi_1^2(16(1-\xi_2^2)^3+9(1-\xi_2^2) -24(1-\xi_2^2)^2)-\\
		& (1-\xi_1^2)(3\xi_2-4\xi_2^2)^2\bigg]\big[\xi_1^2(1-\xi_2^2)-(1-\xi_1^2)\xi_2^2\big]x_1x_2x_3x_4.\end{split}
	\end{align} and finally $$S(\xi_1):=\frac{1}{2}\left(\sqrt[3]{\sqrt{\xi_1^2-1}-\xi_1}+\frac{1}{\sqrt[3]{\sqrt{\xi_1^2-1}-\xi_1}}\right).$$
\end{proposition}
\begin{proof}
	The proof goes analoguous to the proof of Proposition \ref{prop:transform_Mathieu_to_abelian_integral_Sp(N)} or the proof of \cite[Lemma 2.7]{Zwart}. By Lemma \ref{lemma:Euler_angles_G2} and Equation (\ref{eq:finite_type_f_G2}) we have
	\begin{align*}
		\int_{G_2}f(g)^Pdg = C\int_{F^{-1}(G_2)}&\left[\sum_{j=1}^M\sum_{i=1}^Q c_{ij}e^{ik_{ij}^1\widetilde{\phi}_1}\sin^{l_{ij}^1}(\widetilde{\psi}_1)\cos^{m_{ij}^1}(\widetilde{\psi}_1) e^{ik_{ij}^2\widetilde{\omega}_1}e^{ik_{ij}^3\widetilde{\phi}_2}\sin^{l_{ij}^3}(\widetilde{\psi}_2)\cdot\right.\\
		&\cos^{m_{ij}^3}(\widetilde{\psi}_2)e^{ik_{ij}^4\widetilde{\omega}_2}\sin^{l_{ij}^4}(y_1)\cos^{m_{ij}^4}(y_1)\sin^{l_{ij}^5}(y_2)\cos^{m_{ij}^5}(y_2)e^{ik_{ij}^5\phi_1}\cdot\\
		&\left.\sin^{l_{ij}^6}(\psi_1)\cos^{m_{ij}^6}(\psi_1)e^{ik_{ij}^6\omega_1}e^{ik_{ij}^7\phi_2}\sin^{l_{ij}^{7}}(\psi_2)\cos^{m_{ij}^{7}}(\psi_2)e^{i k_{ij}^8\omega_2}\right]^P \cdot\\
		&16 J_{G_2}(y_1,y_2)\sin(\tilde{\psi}_1)\cos(\tilde{\psi}_1)\sin(\tilde{\psi}_2)\cos(\tilde{\psi}_2)\sin(\psi_1)\cos(\psi_1)\cdot\\
		&\sin(\psi_2)\cos(\psi_2)d\tilde{\phi}_1\ldots d\omega_2
	\end{align*}
	We will follow the steps as in \cite{Zwart}. In short, we use the multinomium theorem to expand the exponent of the sum of matrix components, and then we use the following equations $$\int_{0}^{2\pi}e^{i\frac{k}{l}x}dx = \frac{1}{i}\int_{S^*}z^{\frac{k}{l}}\frac{dz}{z},\qquad\qquad\int_{0}^{\pi/2}\sin^{k+p}(\phi)\cos^{l+q}(\phi)d\phi = \int_{0}^1 x^{k+p} (1-x^2)^{\frac{l+q-1}{2}}dx$$ for $p,q,k,l\in\N_0$ with $l>0$ to get to integrals over $S^*$ and $[0,1]$. Then we use the multinomial of Newton again to get the desired result. Note that there are a few integrals of the form $\int_{0}^\pi e^{ikx}dx$, one integral of the form $\int_{0}^{\pi/2} e^{ikx}dx$ and one of the form $\int_{0}^{\pi/4} \cos^{k}(x)\sin^l(x)dx$. The first integral can be transformed to an integral over $S^*$ by subsituting $u=2x$, which gives
	\begin{align*}
		\int_{0}^\pi e^{ikx}dx = \frac{1}{2}\int_{0}^{2\pi} e^{i\frac{k}{2}u}du = \frac{1}{2i}\int_{S^*}z^{\frac{k}{2}}\frac{dz}{z}.
	\end{align*}
	In the same way, we see that a subsitution of $v=4x$ gives
	\begin{align*}
	\int_{0}^{\pi/2} e^{ikx}dx = \frac{1}{4}\int_{0}^{2\pi} e^{i\frac{k}{4}v}dv = \frac{1}{4i}\int_{S^*}z^{\frac{k}{4}}\frac{dz}{z}.
	\end{align*}
	Finally, the last integral can be transformed to 
	\begin{align*}
		\int_{0}^{\pi/4} \cos^{k}(x)\sin^l(x) dx = \int_0^{\frac{1}{\sqrt{2}}}x^k(1-x^2)^{\frac{l-1}{2}}dx=\frac{1}{\sqrt{2}}\int_{0}^1 \left(\frac{x}{\sqrt{2}}\right)^{k}\left(1-\frac{x^2}{2}\right)^{\frac{l-1}{2}}dx.
	\end{align*}
	This allows us to transform all the integrals, which give
	\begin{align*}
		\int_{G_2}f(g)^Pdg &= \frac{C}{\sqrt{2}\cdot2^{7}}\int_{(S^*)^{8}}\int_{[0,1]^4}\int_{\mathcal{A}}\left[\sum_{j=1}^M\sum_{i=1}^Q c_{ij}z_1^{\frac{k_{ij}^1}{2}}\left(\frac{x_1}{\sqrt{2}}\right)^{l_{ij}^1}\left(1-\frac{x_1^2}{2}\right)^{\frac{m_{ij}^1}{2}} z_2^{\frac{k_{ij}^2}{4}}z_3^{\frac{k_{ij}^3}{2}}x_2^{l_{ij}^3}\cdot\right.\\
		&\qquad (1-x_2^2)^{\frac{m_{ij}^3}{2}}z_4^{k_{ij}^4}\sin^{l_{ij}^4}(y_1)\cos^{m_{ij}^4}(y_1)\sin^{l_{ij}^5}(y_2)\cos^{m_{ij}^5}(y_2)z_5^{\frac{k_{ij}^5}{2}}x_3^{l_{ij}^6}(1-x_3^2)^{\frac{m_{ij}^6}{2}}\cdot\\
		&\qquad \left.z_6^{\frac{k_{ij}^6}{2}}z_7^{\frac{k_{ij}^{7}}{2}}x_4^{l_{ij}^{7}}(1-x_4^2)^{\frac{m_{ij}^{7}}{2}}z_8^{k_{ij}^8}\right]^PJ_{G_2}(y_1,y_2)x_1x_2x_3x_4\,dy_1\,dy_2 \,dx_1\ldots dx_4 \frac{dz_1}{z}\ldots \frac{dz_8}{z_8}.
	\end{align*}
	Finally we would like to change the integrals going over $y_1$ and $y_2$. Putting in the area we found for $\mathcal{A}$, we get that
	\begin{align*}
		\int_{G_2}f(g)^Pdg &= \frac{C}{\sqrt{2}\cdot2^{7}}\int_{(S^*)^{8}}\int_{[0,1]^4}\int_{0}^{\pi/2}\int_0^{y_1/3}\left[\sum_{j=1}^M\sum_{i=1}^Q c_{ij}z_1^{\frac{k_{ij}^1}{2}}\left(\frac{x_1}{\sqrt{2}}\right)^{l_{ij}^1}\left(1-\frac{x_1^2}{2}\right)^{\frac{m_{ij}^1}{2}} z_2^{\frac{k_{ij}^2}{4}}z_3^{\frac{k_{ij}^3}{2}}\cdot\right.\\
		&\qquad x_2^{l_{ij}^3}(1-x_2^2)^{\frac{m_{ij}^3}{2}}z_4^{k_{ij}^4}\sin^{l_{ij}^4}(y_1)\cos^{m_{ij}^4}(y_1)\sin^{l_{ij}^5}(y_2)\cos^{m_{ij}^5}(y_2)z_5^{\frac{k_{ij}^5}{2}}x_3^{l_{ij}^6}(1-x_3^2)^{\frac{m_{ij}^6}{2}}\cdot\\
		&\qquad \left.z_6^{\frac{k_{ij}^6}{2}}z_7^{\frac{k_{ij}^{7}}{2}}x_4^{l_{ij}^{7}}(1-x_4^2)^{\frac{m_{ij}^{7}}{2}}z_8^{k_{ij}^8}\right]^PJ_{G_2}(y_1,y_2)x_1x_2x_3x_4\,dy_2\,dy_1 \,dx_1\ldots dx_4 \frac{dz_1}{z}\ldots \frac{dz_8}{z_8}\\
		&= \frac{C}{\sqrt{2}\cdot2^{7}}\int_{(S^*)^{8}}\int_{[0,1]^4}\int_{0}^{\pi/2}\int_0^{\sin(y_1/3)}\left[\sum_{j=1}^M\sum_{i=1}^Q c_{ij}z_1^{\frac{k_{ij}^1}{2}}\left(\frac{x_1}{\sqrt{2}}\right)^{l_{ij}^1}\left(1-\frac{x_1^2}{2}\right)^{\frac{m_{ij}^1}{2}} z_2^{\frac{k_{ij}^2}{4}}z_3^{\frac{k_{ij}^3}{2}}\cdot\right.\\
		&\qquad x_2^{l_{ij}^3}(1-x_2^2)^{\frac{m_{ij}^3}{2}}z_4^{k_{ij}^4}\sin^{l_{ij}^4}(y_1)\cos^{m_{ij}^4}(y_1)\xi_2^{l_{ij}^5}(1-\xi_2^2)^{\frac{m_{ij}^5}{2}}z_5^{\frac{k_{ij}^5}{2}}x_3^{l_{ij}^6}(1-x_3^2)^{\frac{m_{ij}^6}{2}}z_6^{\frac{k_{ij}^6}{2}}z_7^{\frac{k_{ij}^{7}}{2}}\cdot\\
		&\qquad \left.x_4^{l_{ij}^{7}}(1-x_4^2)^{\frac{m_{ij}^{7}}{2}}z_8^{k_{ij}^8}\right]^Pf(y_1,\xi_2)x_1x_2x_3x_4d\xi_2\,dy_1 \,dx_1\ldots dx_4 \frac{dz_1}{z}\ldots \frac{dz_8}{z_8}.
	\end{align*}
	where
	\begin{align*}
		f(y_1,\xi_2)=& 2 \xi_2\sin(2y_1)\left[\sin^2(y_1)(16(1-\xi2^2)^3+9(1-\xi_2^2) -24(1-\xi^2)^2) - \cos^2(y_1)(3\xi_2-4\xi_2^2)^2\right]\cdot\\
		&[\sin^2(y_1)(1-\xi_2^2)-\xi_2^2\cos^2(y_1)]
	\end{align*}
	Note now that doing the subsitution $\xi_1=\sin(y_1)$ is possible, with only the boundary of the $\xi_2$ integral is not easily transformable. However, it follows that $$\sin(y_1/3)=\frac{1}{2}\left(\sqrt[3]{\sqrt{\sin^2(y_1)-1}-\sin(y_1)}+\frac{1}{\sqrt[3]{\sqrt{\sin^2(y_1)-1}-\sin(y_1)}}\right).$$ Let us call the right-handed side $S(\sin(y_1))$, then we see that $\sin(y_1/3)=S(\sin(y_1))$ which is some function depending on $\sin(x)$. Therefore we can apply the subsitution, and we this get the final form
	\begin{align*}
		\int_{G_2}f(g)^Pdg &= \frac{C}{\sqrt{2}\cdot2^{7}}\int_{(S^*)^{8}}\int_{[0,1]^5}\int_0^{S(\xi_1)}\left[\sum_{j=1}^M\sum_{i=1}^Q c_{ij}z_1^{\frac{k_{ij}^1}{2}}\left(\frac{x_1}{\sqrt{2}}\right)^{l_{ij}^1}\left(1-\frac{x_1^2}{2}\right)^{\frac{m_{ij}^1}{2}} z_2^{\frac{k_{ij}^2}{4}}z_3^{\frac{k_{ij}^3}{2}}x_2^{l_{ij}^3}\cdot\right.\\
		&\qquad (1-x_2^2)^{\frac{m_{ij}^3}{2}}z_4^{k_{ij}^4}\xi_1^{l_{ij}^4}(1-\xi_1^2)^{\frac{m_{ij}^4}{2}}\xi_2^{l_{ij}^5}(1-\xi_2^2)^{\frac{m_{ij}^5}{2}}z_5^{\frac{k_{ij}^5}{2}}x_3^{l_{ij}^6}(1-x_3^2)^{\frac{m_{ij}^6}{2}}z_6^{\frac{k_{ij}^6}{2}}z_7^{\frac{k_{ij}^{7}}{2}}x_4^{l_{ij}^{7}}\cdot\\
		&\qquad \left.(1-x_4^2)^{\frac{m_{ij}^{7}}{2}}z_8^{k_{ij}^8}\right]^P\tilde{J}_{G_2}(\xi_1,\xi_2)x_1x_2x_3x_4d\xi_2\,d\xi_1 \,dx_1\ldots dx_4 \frac{dz_1}{z}\ldots \frac{dz_8}{z_8}
	\end{align*}
	where
	\begin{align*}
	\tilde{J}_{G_2}(\xi_1,\xi_2)=& 4 \xi_1\xi_2\left[\xi_1^2(16(1-\xi_2^2)^3+9(1-\xi_2^2) -24(1-\xi_2^2)^2) - (1-\xi_1^2)(3\xi_2-4\xi_2^2)^2\right]\cdot\\
	&[\xi_1^2(1-\xi_2^2)-(1-\xi_1^2)\xi_2^2]
	\end{align*}
	which proves the lemma after noting $\frac{4C}{\sqrt{2}\cdot2^{7}}=\frac{C}{\sqrt{2}\cdot 2^5}=C'$.
\end{proof}
\begin{remark}
	One might wonder why we chose to write the dependencies of $\tilde{\phi}_1,\ldots \omega_2$ in $f$ in the way we did, i.e. why we chose a sum over $e^{ik\tilde{\phi}_1}$ instead of $\sin^{l}(\tilde{\phi}_1)\cos^m(\tilde{\phi}_1)$ when the integrals that come up in Proposition \ref{prop:transform_Mathieu_to_abelian_integral_G_2} are completely different. For in the end, the choice of $e^{ik\tilde{\phi}_1}$ over $\sin^{l}(\tilde{\phi}_1)\cos^m(\tilde{\phi}_1)$ should not change the finite-type function. The idea is to write those parameters as $e^{ikx}$ only if the Jacobian $J$ is independent of $x$, and keep the rest as $\sin^{l}(x)\cos^m(x)$. This way the subsitutions are easily done.
\end{remark}

In the same spirit as in previous section, we see that $\tilde{f}_{G_2}$ is a $\frac{1}{4}$-admissible function. To solve the Mathieu conjecture on $G_2$, we assume the following conjecture:


\begin{conjecture}\label{con:xz-conjecture_G_2}
	Let $f:[0,1]^{6}\times (S^*)^{8}\rightarrow\C$ be a $\frac{1}{4}$-admissible function. If $$\int_{(S^*)^8}\int_{[0,1]^5}\int_0^{S(x_5)}f^P \tilde{J}_{G_2} = 0$$ for all $P\in \N$, then $\vec{0}$ does not lie in the convex hull of $\mathrm{Sp}(f)$.
\end{conjecture}

\begin{theorem}
	Assume Conjecture \ref{con:xz-conjecture_G_2} is true. Then the Mathieu conjecture is true for $G_2$.
\end{theorem}
\begin{proof}
	Let $f,h$ be finite-type functions such that $\int_G f^P(g)dg=0$ for all $P\in\N$. We are interested in $\int_G f^P(g)h(g)dg.$ As discussed, $f$ and $h$ have the form of Equation (\ref{eq:finite_type_f_G2}). By linearity we can again assume $h$ to be of the form
	\begin{align*}
		h(g)=&ce^{iK_1\tilde{\phi}_1}\sin^{L_1}(\tilde{\psi}_1)\cos^{M_1}(\tilde{\psi}_1) e^{iK_2\tilde{\omega}_1}e^{iK_3\tilde{\phi}_2}\sin^{L_3}(\tilde{\psi}_2)\cos^{M_3}(\tilde{\psi}_2)\cdot\\
		&e^{iK_4\tilde{\omega}_2}\sin^{L_4}(y_1)\cos^{M_4}(y_1)\sin^{L_5}(y_2)\cos^{M_5}(y_2)e^{iK_5\phi_1}\sin^{L_6}(\psi_1)\cos^{M_6}(\psi_1)\cdot\\
		&e^{iK_6\omega_1}e^{iK_7\phi_2}\sin^{L_{7}}(\psi_2)\cos^{M_{7}}(\psi_2)e^{i K_8\omega_2}.
	\end{align*}
	Using the multinomial theorem, we get that $\int_{G_2}f^P(g)h(g)dg$ can be written as
	\begin{align}\label{eq:Mathieu_G_2_proof_huge_integral_nonzero}
		\begin{split}\int_{G_2}f^P(g)h(g)dg=& \sum_{\sum \beta_{ij}=P} \binom{P}{\beta_{1,1},\ldots,\beta_{M,Q}}\int_{G_2} c\cdot \left(\prod_{j=1}^M\prod_{i=1}^Qc_{ij}^{\beta_{ij}}\right)e^{i\left(\sum_{i,j}\beta_{ij}k_{ij}^1+K_1\right)\tilde{\phi}_1}\cdot\\
		&\sin^{\sum_{i,j}\beta_{ij}l_{ij}^1+L_1}(\tilde{\psi}_1)\cos^{\sum_{i,j}\beta_{ij}m_{ij}^1+M_1}(\tilde{\psi}_1) e^{i\left(\sum_{i,j}\beta_{ij}k_{ij}^2+K_2\right)\tilde{\omega}_1}e^{i\left(\sum_{i,j}k_{ij}^3+K_3\right)\tilde{\phi}_2}\cdot\\
		&\sin^{\sum_{i,j}\beta_{ij}l_{ij}^3+L_3}(\tilde{\psi}_2)\cos^{\sum_{i,j}\beta_{ij}m_{ij}+M_3}(\tilde{\psi}_2)e^{i\left(\sum_{i,j}\beta_{ij}k_{ij}^4+K_4\right)\tilde{\omega}_2}\sin^{\sum_{i,j}\beta_{ij}l_{ij}^4+L_4}(y_1)\cdot\\
		&\cos^{\sum_{i,j}\beta_{ij}m_{ij}^4+M_4}(y_1)\sin^{\sum_{i,j}\beta_{ij}l_{ij}^5+L_5}(y_2)\cos^{\sum_{i,j}\beta_{ij}m_{ij}^5+M_5}(y_2)e^{i\left(\sum_{i,j}\beta_{ij}k_{ij}^5+K_5\right)\phi_1}\cdot\\
		&\sin^{\sum_{i,j}\beta_{ij}l_{ij}^6+L_6}(\psi_1)\cos^{\sum_{i,j}\beta_{ij}m_{ij}^6+M_6}(\psi_1)e^{i\left(\sum_{i,j}\beta_{ij}k_{ij}^6+K_6\right)\omega_1}e^{i\left(\sum_{i,j}\beta_{ij}k_{ij}^7+K_7\right)\phi_2}\\
		&\sin^{\sum_{i,j}\beta_{ij}l_{ij}^7+L_{7}}(\psi_2)\cos^{\sum_{i,j}\beta_{ij}m_{ij}^7+M_{7}}(\psi_2)e^{i \left(\sum_{i,j}\beta_{ij}k_{ij}^8+K_8\right)\omega_2} J_{G_2}(y_1,y_2) \,d\phi_1\ldots d\tilde{\omega}_2
		\end{split}
	\end{align}
	Let us assume $\int_{G_2}f^P(g)h(g)dg\neq 0$ for infinitely many $P$. Then there exists at least one set of non-negative integers $\{\beta_{ij}\}_{i,j}$ such that $\sum_{i,j}\beta_{ij}=P$ and such that the integral in Equation (\ref{eq:Mathieu_G_2_proof_huge_integral_nonzero}) is non-zero. Applying Proposition \ref{prop:transform_Mathieu_to_abelian_integral_G_2} to this integral, we thus get that
	\begin{align*}
		0\neq \int_{[0,1]^l}\int_{(S^*)^7}& z_1^{\frac{\sum_{i,j}\beta_{ij}k_{ij}^1+K_1}{2}}\left(\frac{x_1}{\sqrt{2}}\right)^{\sum_{i,j}\beta_{ij}l_{ij}^1+L_1}\left(1-\frac{x_1^2}{2}\right)^{\frac{\sum_{i,j}\beta_{ij}m_{ij}^1+M_1}{2}} z_2^{\frac{\sum_{i,j}\beta_{ij}k_{ij}^2+K_2}{4}}z_3^{\frac{\sum_{i,j}k_{ij}^3+K_3}{2}}\\
		&x_2^{\sum_{i,j}\beta_{ij}l_{ij}^3+L_3}(1-x_2^2)^{\frac{\sum_{i,j}\beta_{ij}m_{ij}+M_3}{2}}z_4^{\sum_{i,j}\beta_{ij}k_{ij}^4+K_4}\xi_1^{\sum_{i,j}\beta_{ij}l_{ij}^4+L_4}\cdot\\
		&(1-\xi_1^2)^{\frac{\sum_{i,j}\beta_{ij}m_{ij}^4+M_4}{2}}\xi_2^{\sum_{i,j}\beta_{ij}l_{ij}^5+L_5}(1-\xi_2^2)^{\frac{\sum_{i,j}\beta_{ij}m_{ij}^5+M_5}{2}}z_5^{\frac{\sum_{i,j}\beta_{ij}k_{ij}^5+K_5}{2}}\cdot\\
		&x_3^{\sum_{i,j}\beta_{ij}l_{ij}^6+L_6}(1-x_3^2)^{\frac{\sum_{i,j}\beta_{ij}m_{ij}^6+M_6}{2}}z_6^{\frac{\sum_{i,j}\beta_{ij}k_{ij}^6+K_6}{2}}z_7^{\frac{\sum_{i,j}\beta_{ij}k_{ij}^7+K_7}{2}}\\
		&x_4^{\sum_{i,j}\beta_{ij}l_{ij}^7+L_{7}}(1-x_4^2)^{\frac{\sum_{i,j}\beta_{ij}m_{ij}^7+M_{7}}{2}}z_8^{\sum_{i,j}\beta_{ij}k_{ij}^8+K_8} \tilde{J}_{G_2}(x_1,\ldots,\xi_2) d\xi_1\ldots \frac{dz_8}{z_8}.
	\end{align*}
	To conclude the proof, we will show that $\sum_{i,j}\beta_{ij}k_{ij}^a+K_a=0$ for all $a=1,\ldots,8$ and then make use of Conjecture \ref{con:xz-conjecture_G_2}. To do this, we will make extensive use of the Haar measure properties, i.e. make use of $$\int_{G_2}f(g)dg=\int_{G_2}f(gk)dg=\int_{G_2}f(kg)dg$$ for any measurable function $f$ and any $k\in G_2$. We will denote this as that the functions $g\mapsto gk$ and the function $g\mapsto kg$ are invariant mappings. We make extensive use of our parametrization given in Lemma \ref{lemma:Euler_angles_G2} to see what the invariances of these functions mean for the possible parameters of Equation (\ref{eq:Mathieu_G_2_proof_huge_integral_nonzero}).
	
	We note that sending $g\mapsto \exp(t\lambda_3)g$ is an invariant mapping for all $t\in\R$. Going over to our parametrization this is equivalent to $F_{G_2}(\tilde{\phi}_1,\tilde{\phi}_2,\ldots,\omega_2)\mapsto F_{G_2}(\tilde{\phi}_1+t,\tilde{\phi}_2,\ldots,\omega_2)$ by bijectivity of $F_{G_2}$ up to a measure zero set, which would be equivalent to stating $\tilde{\phi}_1\mapsto \tilde{\phi}_1+t$ is an invariant mapping for the integral given in Equation (\ref{eq:Mathieu_G_2_proof_huge_integral_nonzero}). Since $t$ was chosen arbitrarily, this can only be true if $$\sum_{i,j}\beta_{ij}k_{ij}^1+K_1=0.$$
	The same argument holds that $g\mapsto g \exp(t\lambda_8)$ is an invariant mapping. By the same arguments this is equivalent to $\omega_2\mapsto \omega_2+t$ is an invariant mapping which again can only be true if
	$$\sum_{i,j}\beta_{ij}k_{ij}^8+K_8=0.$$
	Next, we remember that
	\begin{align*}
		dg = J(x,y)\,dg_{K/M}\, dy_1\,dy_2 dk
	\end{align*} where $dy_1dy_2$ is the measure on $\mathcal{A}$, $dg_{K/M}$ is the unique left-invariant measure on $K/M$ and $dk$ is the Haar measure on $K$. This means that if we parametrize $g\in G_2$ as $g=xak$ where $x\in K/M$, $a=\exp(H)\in \exp(\mathcal{A})$ and $k\in K$, then for any measurable function $f$ we have
	\begin{align*}
		\int_{G_2}f(g)dg = \int_{K/M}\int_{A}\left(\int_{K}f(xak)dk\right)dy_1dy_2\,dg_{K/M} &= \int_{K/M}\int_{A}\left(\int_{K}f(xahk)dk\right)dy_1dy_2\,dg_{K/M}\\
		&=\int_{K/M}\int_{A}\left(\int_{K}f(xakh)dk\right)dy_1dy_2\,dg_{K/M}.
	\end{align*}
	for any $h\in K$. In other words, the mapping $g=xak\mapsto xahk$ for any $h\in K$ is also an invariant mapping. This way, the mapping $g=xak\mapsto xa\exp(t\lambda_{3})k$ is equivalent to sending $\phi_1\mapsto \phi_1+t$ which is invariant. This can only be the case in Equation (\ref{eq:Mathieu_G_2_proof_huge_integral_nonzero}) if $$\sum_{i,j}\beta_{ij}k_{ij}^5+K_5=0.$$
	Next, remember that $K \simeq SU(2)\times SU(2)/\mathbb{Z}_2$. In other words, it is the product of two commuting groups. So $$e^{t\lambda_8}k=e^{\phi_1\lambda_3}e^{\psi_1\lambda_2}e^{\omega_1\lambda_3}e^{t\lambda_8}e^{\phi_2\lambda_8}e^{\psi_2\lambda_9}e^{\omega_2\lambda_8}=e^{\phi_1\lambda_3}e^{\psi_1\lambda_2}e^{\omega_1\lambda_3}e^{(\phi_2+t)\lambda_8}e^{\psi_2\lambda_9}e^{\omega_2\lambda_8}$$ if $k=e^{\phi_1\lambda_3}e^{\psi_1\lambda_2}e^{\omega_1\lambda_3}e^{\phi_2\lambda_8}e^{\psi_2\lambda_9}e^{\omega_2\lambda_8}$. Thus by previous arguments we see that the mapping $g=xak\mapsto xa\exp(t\lambda_{8})k$ is equivalent to the map $\phi_2\mapsto \phi_2+t$ being invariant, which can only be the case if $$\sum_{i,j}\beta_{ij}k_{ij}^7+K_7=0.$$ 
	Finally, we note that the measure $dg_{K/M}=dk$ as well, hence the same arguments can be reused to argue that the mapping $g=xak\mapsto xhak$ for any $h\in K$ is also an invariant mapping. Applying that to $g=xak\mapsto x\exp(t\lambda_3)ak$ is then equivalent to sending $\tilde{\omega}_1\mapsto \tilde{\omega}_1+t$ which can only be the case if $$\sum_{i,j}\beta_{ij}k_{ij}^2+K_2=0.$$ 	In the same way, the mapping $g=xak\mapsto x\exp(t\lambda_{8})ak$ is equivalent to sending $\tilde{\omega}_2\mapsto \tilde{\omega}_2+t$ which should be invariant, meaning that $$\sum_{i,j}\beta_{ij}k_{ij}^4+K_4=0,$$
	and the map $g\mapsto g\exp(t\lambda_3)$ is invariant which is equivalent to sending $\tilde{\omega}_1\mapsto \tilde{\omega}_1+t$ which can only be true if $$\sum_{i,j}\beta_{ij}k_{ij}^6+K_6=0.$$
	Finally, using the invariance of the mapping $g\mapsto \exp(t\lambda_8)g$ being equivalent to $\tilde{\phi}_2\mapsto \tilde{\phi}_2+t$ being invariant, we can conclude that $$\sum_{i,j}\beta_{ij}k_{ij}^3+K_3=0.$$
	In other words, we have that $\sum_{i,j}\beta_{ij}k_{ij}^a+K_a=0$ for all $a=1,\ldots 8$. This means that $$-(K_1,\ldots,K_8)=\sum_{i,j}\beta_{ij}(k_{ij}^1,\ldots,k_{ij}^8).$$ Diving both terms by $P\in\N$ gives $$\frac{1}{P}(K_1,\ldots,K_8)=\frac{1}{P}\sum_{i,j}\beta_{ij}(k_{ij}^1,\ldots,k_{ij}^8).$$ Taking the limit $P\rightarrow\infty$ gives that the left hand side goes to 0. Hence there exists a limiting sequence in $\mathrm{Sp}(f)$ that converges to 0. Since the convex hull is a closed set, it means $0\in \mathrm{Conv}(\mathrm{Sp}(f))$ which is a contradiction with out assumption. So it cannot be that $\int_{G_2}f^Ph\neq 0$ for infintely many $P$, hence proving the theorem.
	
\end{proof}

\newpage

\appendix
\section[Proof of Theorem \ref{thm:KAK_decomp}]{Proof of Theorem \ref{thm:KAK_decomp}}\label{appendix:Proof_General_KAK_decomp}
In this section, we prove Theorem \ref{thm:KAK_decomp}. For completeness we restate the theorem. The proof is based on \cite{Helgason_DifGeom_Lie_Groups_Symm_Spaces,KnappBeyond}.

\KAKDecomp*

\begin{proof}
	We define the mapping $$f:K/M\times \exp(\mathcal{A})\times K\rightarrow G,\qquad f(kM,\exp(H),l):=\exp( \Ad(k)H)l.$$ We claim that this mapping is a surjection, and up to a measure zero set a diffeomorphism.
	
	First we note that $f$ is well-defined. Indeed, taking any two representatives $k,l\in kM$ gives $k=lm$ with $m\in M$. Then it is clear $\Ad(k)H=\Ad(l)\Ad(m)H=\Ad(l)H$. In addition, we note that $f(kM,\exp(H),l)=k\exp(H)k^{-1}l$. According to Helgason \cite[Thm. VII.8.6]{Helgason_DifGeom_Lie_Groups_Symm_Spaces} $$G=K\exp(\mathcal{A})K,$$ and if $k\exp(H)l\in K\exp(\mathcal{A})K$ then $H$ is unique. This show surjectivity of $f$.
	Now assume that $$u:=f(k_1M,\exp(H_1),l_1)=f(k_2M,\exp(H_2),l_2).$$ Then $H_1=H_2=:H$ as mentioned before, and so we are left with 
	\begin{align*}
		\exp(\Ad(k_1)H)l_1 &=\exp(\Ad(k_2)H)l_2\\
		k_1 \exp(H)k_1^{-1}l_1&= k_2\exp(H)k_2^{-1}l_2.
	\end{align*}
	We remember that $\theta:\mathfrak{g}\rightarrow \mathfrak{g}$ is an automorphism, and $G$ is simply-connected, hence we can lift $\theta$ to an automorphism $\Theta:G\rightarrow G$ by $\Theta(\exp(X))=\exp(\theta X)$. Then $\Theta(a)=a^{-1}$ for all $a\in \exp\mathfrak{p}$ and $\Theta(k)=k$ for all $k\in K$. Therefore we see that 
	\begin{align*}
		u(\Theta(u))^{-1}&=k_1\exp(H)k_1^{-1}l_1l_1^{-1}k_1\exp(H)k_1^{-1} = k_1\exp(2H)k_1^{-1}\\
		&= k_2\exp(2H)k_2^{-1}.
	\end{align*}
	The last equality sign shows that
	\begin{align*}
		\exp(2H)&=k_2^{-1}k_1\exp(2H)[k_2^{-1}k_1]^{-1}\\
		&= \exp(\Ad(k)2H)
	\end{align*}
	where we defined $k:=k_2^{-1}k_1$. 
	
	Let us now assume that $\alpha(H)\notin \pi i\Z$ for all $\alpha\in \Delta_\mathfrak{p}$, i.e. $H\in\mathrm{int}(\mathcal{A})$. Using a similar argument as in Knapp \cite[Thm. 7.36]{KnappBeyond} we get that the equality $\exp(\Ad(k)2H)=\exp(2H)$ requires $\Ad(k)H=H$. 
	
	We claim that $Z_\mathfrak{g}(H)=\mathfrak{m}\oplus \mathfrak{a}$ where $\mathfrak{m}=Z_\mathfrak{k}(\mathfrak{a})$. For if $X\in Z_\mathfrak{g}(H)$, we see that $X\in \mathfrak{g}_\C$ can be written as $$X=X_0+\sum_{\alpha\in\Delta_{\mathfrak{p}}}X_\alpha$$ where $X_\alpha\in\mathfrak{g}_\alpha$ and $X_0\in\mathfrak{g}_0=\mathfrak{m}_\C\oplus\mathfrak{a}_\C$. Then we note that $$\ad(H)X=0+\sum_{\alpha\in\Delta_\mathfrak{p}}\alpha(H)X_\alpha$$ but $\ad(H)X=0$ by definition of $X$. So $\alpha(H)X_\alpha=0$. However $H\in\mathfrak{a}_+$ so $X_\alpha = 0$. This shows that $X\in \mathfrak{m}\oplus\mathfrak{a}$ and so $Z_\mathfrak{g}(H)=\mathfrak{m}\oplus \mathfrak{a}$. In particular $Z_\mathfrak{p}(H)=\mathfrak{a}$. In this particular case it means that $$\Ad(k)\mathfrak{a} = Z_\mathfrak{p}(\Ad(k)H) = Z_\mathfrak{p}(H) = \mathfrak{a}.$$ So $\Ad(k)\mathfrak{a}=\mathfrak{a}$ and so $k\in N_K(\mathfrak{a})$.
	
	Since $H\in\mathfrak{a}_+$, using \cite[Lemma VII.2.2]{Helgason_DifGeom_Symm_Spaces} there exists an $s\in W(U,K)$ such that $\Ad(k)H=s.H$, where $W(U,K):=N_K(\mathfrak{a})/Z_K(\mathfrak{a})$ is the analytic Weyl group. Since $H$ lies in the positive Weyl chamber, the only reflection satisfying $s.H=H$ is $s=1$. So $\Ad(k)=1$, which means that $k\in Z_K(\mathfrak{a})$. 
	
	Going back to the definition of $k$, we thus see that $k_2=k_1k$ with $k\in Z_K(\mathfrak{a})=M$. In other words $k_2M=k_1M$. It is them immediate that $l_1=l_2$ and therefore we find that $f$ is a bijection on $\mathrm{int}(\mathcal{A})$. Note that $\mathcal{A}$ and $\mathrm{int}(\mathcal{A})$ differ only by a measure zero set.
	
	Note that $f$ is smooth. This can be seen by first noting that the map $\psi:K\times \mathcal{A}\rightarrow \mathfrak{g}$ given by $(k,H)\mapsto \Ad(k)H$ is smooth. Hence by the universal properties of the quotient, the map $\tilde{\psi}:K/M\times \mathcal{A}\rightarrow \mathfrak{g}$, given by $(kM,H)\mapsto \Ad(k)H$, is also smooth. Therefore we see that the map $\exp\circ \tilde{\psi}$ is a smooth map, and multiplication is smooth, hence $f$ is smooth. 
	
	So $f: K/M \times \exp(\mathrm{int}(\mathcal{A}))\times K\rightarrow G$ is bijective and smooth. We only need to show that the inverse is smooth, or that the tangent map is bijective everywhere, i.e. $f$ is regular. We will show the latter. Let $\pi: K\rightarrow K/M$ and $\rho: G\rightarrow G/K$ the canonical maps. We remind ourselves that the map $$\Phi: K/M\times \mathcal{A}\rightarrow G/K,\qquad (kM,H)\mapsto \rho[\exp(\Ad_G(k)H)]$$ is a regular map at every point $(kM,H)\in K/M\times \mathrm{int}(\mathcal{A})$ \cite[Prop. VII.3.2]{Helgason_DifGeom_Lie_Groups_Symm_Spaces}. Now note that on $\mathrm{int}(\mathcal{A})$ the exponential $\exp:\mathrm{int}(\mathcal{A})\rightarrow \exp(\mathrm{int}(\mathcal{A}))$ is a diffeomorphism. To prove this, note that since $\mathfrak{a}$ is abelian, the differential of the exponent is just the identity. It is also injective, for if $\exp(H)=\exp(H')$ then $\exp(H-H')=e$. This means $\Ad(\exp(H-H'))$ is the identity mapping on $\mathfrak{g}$ which in turn means $\alpha(H)-\alpha(H')\in 2\pi i\Z$ for all $\alpha\in \Delta_\mathfrak{p}$. If there exists an $\alpha\in \Delta_\mathfrak{p}$ such that $\alpha(H)-\alpha(H') = 2\pi i n$ with $n\in\mathbb{Z}-\{0\}$ then note that either $H$ or $H'$ is not in $\mathcal{A}$ since it means $|\alpha(H)|$ or $|\alpha(H')|$ is bigger or equal than $\pi$ which is impossible by definition of $\mathrm{int}(\mathcal{A})$. If $\alpha(H)-\alpha(H')=0$ for all $\alpha\in\Delta_\mathfrak{p}$ then note that $\Delta_\mathfrak{p}$ is a root system, this means that $\Delta_\mathfrak{p}$ span $\mathfrak{a}^*$, hence we must have $H-H'=0$. So $H=H'$. Hence we conclude that $$\Phi: K/M\times \exp(\mathrm{int}(\mathcal{A}))\rightarrow G/K,\qquad (kM,\exp(H))\mapsto \rho[\exp(\Ad_G(k)H)]$$ 
	is a regular map at every point $(kM,\exp(H))\in K/M\times \exp(\mathrm{int}(\mathcal{A}))$.
	
	Note that $\Ad(k)H\in\mathfrak{p}$ for all $k\in K$, and $T_{eK}(G/K)\simeq \mathfrak{g}/\mathfrak{k}\simeq \mathfrak{p}$ and thus $T_{gK}(G/K)=\left(T_{eK}\tau(g)\right)\left[ T_{eK}(G/K)\right]$ where $\tau(g): xK\mapsto gxK$ is a diffeomorphism for all $g\in G$. This means that $T_{gK}(G/K)\simeq \mathfrak{p}$ for all $gK\in G/K$. Since $T_{e}\rho: \mathfrak{g}\mapsto \mathfrak{g}/\mathfrak{k}\simeq \mathfrak{p}$ by the same inclusion, we can conclude from this that the mapping $$\tilde{\Phi}:K/M\times \exp(\mathcal{A})\rightarrow \exp(\mathfrak{p}),\qquad (kM,\exp(H))\mapsto \exp(\Ad_G(k)H)$$ is a regular map. Now we note that $$f:K/M\times \exp{\mathcal{A}}\times K\rightarrow G,\qquad (kM,\exp(H),l)\mapsto \exp(\Ad_G(k)H)l$$ can be written as $$f(kM,\exp(H),l)=\mu(\tilde{\Phi}(kM, \exp(H)), l)$$ where $\mu:G\times G\rightarrow G, (g,h)\mapsto gh$ is the multiplication map. In this way, we see that $$T_{kM,\exp(H),l}f(X,Q,Y) = T_{kM,\exp(H)}\tilde{\Phi}(X,Q)+Y$$ which is clearly a bijection, since $\mathfrak{g}=\mathfrak{k}\oplus\mathfrak{p}$. In other words, we have that $f$ is a smooth regular bijection, i.e. a diffeomorphism.
	
	To get the decomposition of the Haar measure, we apply the techniques of \cite[Section X.1.5]{Helgason_DifGeom_Symm_Spaces} which gives that $$\int_G h(g)dg= C \int_{K/M}\int_{\mathcal{A}}\int_K(h\circ f)(k_1M, \exp(H),k_2)\left|\prod_{\alpha\in \Delta_{\mathfrak{p}^+}}\sin(\alpha(iH))\right|dk_2dHdk_M$$ where $C>0$ is some constant. Here $dk_2$ is the Haar measure on $K$, $dH$ is the Haar measure on $\mathfrak{a}$ and $dk_M$ is the unique left-invariant measure on $K/M$. We can rewrite this as 
	\begin{align*}
		\int_G h(g)dg&= C \int_{K/M}\int_{\mathcal{A}}\int_K(h\circ f)(k_1M, \exp(H),k_2)\left|\prod_{\alpha\in \Delta_{\mathfrak{p}^+}}\sin(\alpha(iH))\right|dk_2dHdk_M\\
		&=C \int_{K/M}\int_{\mathcal{A}}\int_Kh\left(\exp(\Ad(k_1)H)k_2\right)\left|\prod_{\alpha\in \Delta_{\mathfrak{p}^+}}\sin(\alpha(iH))\right|dk_2dHdk_M\\
		&=C \int_{K/M}\int_{\mathcal{A}}\int_Kh(k_1\exp(H)k_1^{-1}k_2)\left|\prod_{\alpha\in \Delta_{\mathfrak{p}^+}}\sin(\alpha(iH))\right|dk_2dHdk_M.
	\end{align*}
	Using the fact that $dk_2$ is a Haar measure, we can write this as 
	\begin{align*}
		\int_G h(g)dg&=C \int_{K/M}\int_{\mathcal{A}}\int_Kh(k_1\exp(H)k_2)\left|\prod_{\alpha\in \Delta_{\mathfrak{p}^+}}\sin(\alpha(iH))\right|dk_2dHdk_M.
	\end{align*}
	which proves the theorem
\end{proof}

\section{Some results regarding $G_2$}\label{appendix:spanning_set_G_2}
This section is based on \cite{Cacciatori2}. We can see $G_2$ as the group of automorphisms on the octonions. This already shows that $G_2$ can be embedded into $GL(7,\C)$, but we want to be explicit enough to make sure that the construction done in the main document can be verified explicitly. The basis elements for $\mathfrak{g}_2:=T_eG_2$ are given by
\begin{align*}
	\lambda_1&:=\scalemath{0.8}{\begin{pmatrix}
	0&0&0&0&0&0&0\\
	0&0&0&0&0&0&0\\
	0&0&0&0&0&0&0\\
	0&0&0&0&0&0&-1\\
	0&0&0&0&0&-1&0\\
	0&0&0&0&1&0&0\\
	0&0&0&1&0&0&0
	\end{pmatrix}},\qquad \lambda_2:=\scalemath{0.8}{\begin{pmatrix}
	0&0&0&0&0&0&0\\
	0&0&0&0&0&0&0\\
	0&0&0&0&0&0&0\\
	0&0&0&0&0&1&0\\
	0&0&0&0&0&0&-1\\
	0&0&0&-1&0&0&0\\
	0&0&0&0&1&0&0
	\end{pmatrix}},\qquad\lambda_3:=\scalemath{0.8}{\begin{pmatrix}
	0&0&0&0&0&0&0\\
	0&0&0&0&0&0&0\\
	0&0&0&0&0&0&0\\
	0&0&0&0&-1&0&0\\
	0&0&0&1&0&0&0\\
	0&0&0&0&0&0&-1\\
	0&0&0&0&0&1&0
	\end{pmatrix}},\\
	\lambda_4&:=\scalemath{0.8}{\begin{pmatrix}
	0&0&0&0&0&0&0\\
	0&0&0&0&0&0&1\\
	0&0&0&0&0&1&0\\
	0&0&0&0&0&0&0\\
	0&0&0&0&0&0&0\\
	0&0&-1&0&0&0&0\\
	0&-1&0&0&0&0&0
	\end{pmatrix}},\qquad \lambda_5:=\scalemath{0.8}{\begin{pmatrix}
	0&0&0&0&0&0&0\\
	0&0&0&0&0&-1&0\\
	0&0&0&0&0&0&1\\
	0&0&0&0&0&0&0\\
	0&0&0&0&0&0&0\\
	0&1&0&0&0&0&0\\
	0&0&-1&0&0&0&0
	\end{pmatrix}},\qquad\lambda_6:=\scalemath{0.8}{\begin{pmatrix}
	0&0&0&0&0&0&0\\
	0&0&0&0&1&0&0\\
	0&0&0&-1&0&0&0\\
	0&0&1&0&0&0&0\\
	0&-1&0&0&0&0&0\\
	0&0&0&0&0&0&0\\
	0&0&0&0&0&0&0
	\end{pmatrix}},\\
\lambda_7&:=\scalemath{0.8}{\begin{pmatrix}
	0&0&0&0&0&0&0\\
	0&0&0&-1&0&0&0\\
	0&0&0&0&-1&0&0\\
	0&1&0&0&0&0&0\\
	0&0&1&0&0&0&0\\
	0&0&0&0&0&0&0\\
	0&0&0&0&0&0&0
	\end{pmatrix}},\qquad \lambda_8:=\scalemath{0.8}{\begin{pmatrix}
	0&0&0&0&0&0&0\\
	0&0&-2&0&0&0&0\\
	0&2&0&0&0&0&0\\
	0&0&0&0&1&0&0\\
	0&0&0&-1&0&0&0\\
	0&0&0&0&0&0&-1\\
	0&0&0&0&0&1&0
	\end{pmatrix}},\quad\lambda_9:=\scalemath{0.8}{\begin{pmatrix}
	0&-2&0&0&0&0&0\\
	2&0&0&0&0&0&0\\
	0&0&0&0&0&0&0\\
	0&0&0&0&0&0&1\\
	0&0&0&0&0&-1&0\\
	0&0&0&0&1&0&0\\
	0&0&0&-1&0&0&0
	\end{pmatrix}},\\
\lambda_{10}&:=\scalemath{0.8}{\begin{pmatrix}
	0&0&-2&0&0&0&0\\
	0&0&0&0&0&0&0\\
	2&0&0&0&0&0&0\\
	0&0&0&0&0&-1&0\\
	0&0&0&0&0&0&-1\\
	0&0&0&1&0&0&0\\
	0&0&0&0&1&0&0
	\end{pmatrix}},\quad \lambda_{11}:=\scalemath{0.8}{\begin{pmatrix}
	0&0&0&-2&0&0&0\\
	0&0&0&0&0&0&-1\\
	0&0&0&0&0&1&0\\
	2&0&0&0&0&0&0\\
	0&0&0&0&0&0&0\\
	0&0&-1&0&0&0&0\\
	0&1&0&0&0&0&0
	\end{pmatrix}},\quad\lambda_{12}:=\scalemath{0.8}{\begin{pmatrix}
	0&0&0&0&-2&0&0\\
	0&0&0&0&0&1&0\\
	0&0&0&0&0&0&1\\
	0&0&0&0&0&0&0\\
	2&0&0&0&0&0&0\\
	0&-1&0&0&0&0&0\\
	0&0&-1&0&0&0&0
	\end{pmatrix}},\\
\lambda_{13}&:=\scalemath{0.8}{\begin{pmatrix}
	0&0&0&0&0&-2&0\\
	0&0&0&0&-1&0&0\\
	0&0&0&-1&0&0&0\\
	0&0&1&0&0&0&0\\
	0&1&0&0&0&0&0\\
	2&0&0&0&0&0&0\\
	0&0&0&0&0&0&0
	\end{pmatrix}},\quad \lambda_{14}:=\scalemath{0.8}{\begin{pmatrix}
	0&0&0&0&0&0&-2\\
	0&0&0&1&0&0&0\\
	0&0&0&0&-1&0&0\\
	0&-1&0&0&0&0&0\\
	0&0&1&0&0&0&0\\
	0&0&0&0&0&0&0\\
	2&0&0&0&0&0&0
	\end{pmatrix}}.
\end{align*} 
Next, we briefly discuss the root decomposition of $\mathfrak{g}_2$. It is clear that if we define the involution $$\theta=\Ad\begin{pmatrix}
\mathbf{1}_3&\\
&-\mathbf{1}_4
\end{pmatrix}$$ that $\mathfrak{g}_2=\mathfrak{k}\oplus\mathfrak{p}$ where $$\mathfrak{k}:=\mathrm{span}_\R(\lambda_1,\lambda_2,\lambda_3,\lambda_8,\lambda_9,\lambda_{10})$$ and $$\mathfrak{p}:=\mathrm{span}_\R(\lambda_4,\lambda_5,\lambda_6,\lambda_7,\lambda_{11},\lambda_{12},\lambda_{13},\lambda_{14}).$$ In specific, note that $[\lambda_5,\lambda_{11}]=0$. Choosing $$\mathfrak{a}:=\mathrm{span}_\R(\lambda_5,\lambda_{11})$$ gives a maximal abelian subalgebra in $\mathfrak{p}$. One can compute all commutators, and find that in the basis $\{\lambda_1,\ldots,\lambda_4,\lambda_6,\ldots,\lambda_{10},\lambda_{12},\ldots,\lambda_{14}\}$ the operators $\ad(\lambda_5)$ and $\ad(\lambda_{11})$ look like

\begin{align*}
\setcounter{MaxMatrixCols}{20}
	\ad(\lambda_5)&=\begin{pmatrix}
	0 & 0 & 0 & 0 & 1 & 0 & 0 & 0 & 0 & 0 & 0 & 0\\
	0 & 0 & 0 & 0 & 0 & 1 & 0 & 0 & 0 & 0 & 0 & 0\\
	0 & 0 & 0 & -1 & 0 & 0 & 0 & 0 & 0 & 0 & 0 & 0\\
	0 & 0 & 1 & 0 & 0 & 0 & 3 & 0 & 0 & 0 & 0 & 0\\
	-1 & 0 & 0 & 0 & 0 & 0 & 0 & 0 & 0 & 0 & 0 & 0\\
	0 & -1 & 0 & 0 & 0 & 0 & 0 & 0 & 0 & 0 & 0 & 0\\
	0 & 0 & 0 & -1 & 0 & 0 & 0 & 0 & 0 & 0 & 0 & 0\\
	0 & 0 & 0 & 0 & 0 & 0 & 0 & 0 & 0 & 0 & -1 & 0\\
	0 & 0 & 0 & 0 & 0 & 0 & 0 & 0 & 0 & 0 & 0 & 1\\
	0 & 0 & 0 & 0 & 0 & 0 & 0 & 0 & 0 & 0 & 0 & 0\\
	0 & 0 & 0 & 0 & 0 & 0 & 0 & 1 & 0 & 0 & 0 & 0\\	
	0 & 0 & 0 & 0 & 0 & 0 & 0 & 0 & -1 & 0 & 0 & 0
	\end{pmatrix},\\
	\ad(\lambda_{11})&=\begin{pmatrix}
	0 & 0 & 0 & 0 & 0 & 0 & 0 & 0 & 0 & 0 & 0 & 3\\
	0 & 0 & 0 & 0 & 0 & 0 & 0 & 0 & 0 & 0 & -3 & 0\\
	0 & 0 & 0 & 0 & 0 & 0 & 0 & 0 & 0 & 3 & 0 & 0\\
	0 & 0 & 0 & 0 & 0 & 0 & 0 & 0 & 0 & 0 & 0 & 0\\
	0 & 0 & 0 & 0 & 0 & 0 & 0 & 0 & -3 & 0 & 0 & 0\\
	0 & 0 & 0 & 0 & 0 & 0 & 0 & -3 & 0 & 0 & 0 & 0\\
	0 & 0 & 0 & 0 & 0 & 0 & 0 & 0 & 0 & -1 & 0 & 0\\
	0 & 0 & 0 & 0 & 0 & 1 & 0 & 0 & 0 & 0 & 0 & -2\\
	0 & 0 & 0 & 0 & 1 & 0 & 0 & 0 & 0 & 0 & 2 & 0\\
	0 & 0 & -1 & 0 & 0 & 0 & 1 & 0 & 0 & 0 & 0 & 0\\
	0 & 1 & 0 & 0 & 0 & 0 & 0 & 0 & -2 & 0 & 0 & 0\\
	-1 & 0 & 0 & 0 & 0 & 0 & 0 & 2 & 0 & 0 & 0 & 0
	\end{pmatrix}
\end{align*}
Computing the eigenvalues and eigenvectors of these matrices give the roots of $\mathfrak{g}_2$. Let us define the following linear functionals $\alpha,\beta\in\mathfrak{a}^*$ by
\begin{align*}
	\alpha(\lambda_5)&=0,\qquad \alpha(\lambda_{11})=2\\
	\beta(\lambda_5)&=1,\qquad \beta(\lambda_{11})=-3.
\end{align*}
We claim $i\alpha$ and $i\beta$ are the simple roots for the complex Lie algebra $(\mathfrak{g}_2)_\C$. In fact, the root system is given by $$\Delta_\mathfrak{p}=\{i\alpha,i\beta,\pm i(\alpha+\beta),\pm i(2\alpha+\beta),\pm i(3\alpha+\beta),\pm i(3\alpha+2\beta)
\},$$ as was expected. Choosing the positive roots as all the roots having $+$, we find that $\alpha,\beta$ are the simple roots. The root spaces are given by
\begin{align*}
	\mathfrak{g}_{i\alpha} &= \C(\lambda_3-2i\lambda_4+\lambda_8),\qquad \mathfrak{g}_{-i\alpha}= \C(\lambda_3+2i\lambda_4+\lambda_8),\\
	\mathfrak{g}_{i\beta} &=\C(i\lambda_1 + i\lambda_2 -\lambda_6-\lambda_7-i\lambda_9-\lambda_{10}-\lambda_{13}+\lambda_{14}),\\
	\mathfrak{g}_{i\beta} &=\C(i\lambda_1 - i\lambda_2 -\lambda_6-\lambda_7+i\lambda_9-\lambda_{10}-\lambda_{13}+\lambda_{14}),\\
	\mathfrak{g}_{i(\alpha+\beta)}&= \C(3i\lambda_1-3i\lambda_2-3\lambda_6+3\lambda_7+i\lambda_9-i\lambda_{10}+\lambda_{13}+\lambda_{14}),\\
	\mathfrak{g}_{-i(\alpha+\beta)}&=\C(-3i\lambda_1+3i\lambda_2-3\lambda_6+3\lambda_7-i\lambda_9+i\lambda_{10}+\lambda_{13}+\lambda_{14}),\\
	\mathfrak{g}_{i(2\alpha+\beta)}&= \C(-3i\lambda_1-3i\lambda_2+3\lambda_6+3\lambda_7-i\lambda_9-i\lambda_{10}-\lambda_{12}+\lambda_{14}),\\
	\mathfrak{g}_{-i(2\alpha+\beta)}&=\C(3i\lambda_1+3i\lambda_2+3\lambda_6+3\lambda_7+i\lambda_9+i\lambda_{10}-\lambda_{12}+\lambda_{14}),\\
	\mathfrak{g}_{i(3\alpha+\beta)}&= \C(-i\lambda_1+i\lambda_2+\lambda_6-\lambda_7+i\lambda_9-i\lambda_{10}+\lambda_{13}+\lambda_{14}),\\
	\mathfrak{g}_{-i(3\alpha+\beta)}&=\C(i\lambda_1-i\lambda_2+\lambda_6-\lambda_7-i\lambda_9+i\lambda_{10}+\lambda_{13}+\lambda_{14}),\\
	\mathfrak{g}_{i(3\alpha+2\beta)}&=\C(-3i\lambda_3+i\lambda_8+2\lambda_{12}),\qquad \mathfrak{g}_{-i(3\alpha+2\beta)}=\C(3i\lambda_3-i\lambda_8+2\lambda_{12}).
\end{align*}

\newpage
\bibliography{bibfile}
	
\end{document}